\documentclass[12pt]{amsart}
\usepackage{amsmath}
\usepackage{amsfonts}
\usepackage{amssymb}
\usepackage{graphics}
\usepackage{epstopdf}
\usepackage{epsf}
\usepackage{epsfig}
\usepackage{url}

\newtheorem{theorem}{Theorem}

   \newtheorem{problem}[theorem]{Problem}
   
    \newtheorem{corollary}[theorem]{Corollary}
   \newtheorem{proposition}[theorem]{Proposition}

   \theoremstyle{definition}
   \newtheorem{definition}[theorem]{Definition}
   \newtheorem{example}[theorem]{Example}

\newcommand{\CC}{\mathbb C}
\newcommand{\RR}{\mathbb R}

\newcommand{\PP}{\mathbb P}
\newcommand{\N}{\mathcal N}

\begin{document}

\title[The Hyperdeterminant  and Triangulations of the 4-Cube]{The 
 Hyperdeterminant  and \\ Triangulations of the 4-Cube}

\author[Peter Huggins, Bernd Sturmfels, Josephine Yu,
and Debbie Yuster]{Peter Huggins, 
Bernd Sturmfels, \\ Josephine Yu
and Debbie S. Yuster}

\begin{abstract}
The hyperdeterminant of format $2\times 2 \times 2 \times 2$
is a polynomial of degree $24$ in $16$ unknowns which has
$ 2894276$ terms. We compute the Newton polytope of this
polynomial and the secondary polytope of the $4$-cube.
The $87959448 $  regular triangulations
of the $4$-cube are classified
into $25448$  $D$-equivalence classes,
one for each vertex of the Newton polytope.
The $4$-cube has $80876$ coarsest regular subdivisions,
one for each facet  of the secondary polytope, but
only $268$ of them come from the hyperdeterminant.

\end{abstract}
\maketitle
\section{Introduction}

This article makes three contributions to the mathematics of computation:
\begin{itemize}
\item solution to a computational challenge problem posed by I.M.~Gel'fand,
\item new theorems in algebra (hyperdeterminant) and 
geometry ($4$-cube),
\item new computational methodology 
(algorithms and software) for large polynomials
and large convex polytopes in the presence of symmetry.
\end{itemize}

\smallskip

We begin our discussion with a small familiar polynomial in four unknowns:
$$ D_{22} \quad = \quad c_{00} c_{11} - c_{01} c_{10}. $$  
This irreducible quadratic polynomial is the determinant of a $2 \times 2$-matrix $(c_{ij})$.
Next, the {\em hyperdeterminant} of a $2 \times 2 \times 2$-tensor is the following irreducible
polynomial in eight unknowns which is the sum of
twelve terms of degree four:
\begin{eqnarray*}
D_{222} \quad =
& \quad 4 c_{000}c_{011}c_{101}c_{110}+ 4c_{001}c_{010}c_{100}c_{111} \\
&+ c_{000}^2c_{111}^2 + c_{001}^2c_{110}^2 + c_{010}^2c_{101}^2 +
c_{011}^2c_{100}^2
\\
& -  2 c_{000}c_{001}c_{110}c_{111}
-2 c_{000}c_{010}c_{101}c_{111}
  -2 c_{000}c_{011}c_{100}c_{111} \\
  & -2 c_{001}c_{010}c_{101}c_{110}
   -2 c_{001}c_{011}c_{110}c_{100}
   -2 c_{010}c_{011}c_{101}c_{100}.
     \end{eqnarray*}
This polynomial is  known as the {\em tangle} in the
physics literature. 
Hyperdeterminants of arbitrary format were
studied by Gel'fand, Kapranov, and Zelevinsky
\cite{GKZ}. They  have numerous applications ranging
from quantum information theory \cite{Lev, LT} to
computational biology \cite{AR, SJ} and numerical analysis \cite{DSL, tBe}.
Basic definitions concerning hyperdeterminants will be reviewed in Section 3.

The starting point of this paper is a problem stated by
 I.~M.~Gel'fand in his Fall 2005 research seminar at Rutgers University:
 find the monomial expansion of the
hyperdeterminant $D_{2222}$ of a $2 \times 2 \times 2 \times 2$-tensor.
Back in 1852, Schl\"afli 
had given a nested formula (see \cite[\S 14.4]{GKZ})  for $D_{2222}$, which
 is a polynomial of degree $24$ in  $16$ unknowns
$\,c_{0000}, c_{0001}, \ldots, c_{1111}$, and, 150 years later,
 Luque and Thibon \cite{LT} expressed
$D_{2222}$ in terms of the fundamental tensor invariants.
Gel'fand's challenge to the computer algebra community was to
 expand these formulas into monomials.
We solved Gel'fand's problem:

\begin{theorem}
The hyperdeterminant $D_{2222}$ is the sum of
$ 2894276$ monomials in $9617$ orbits.
The Newton polytope of $D_{2222}$ is an $11$-dimensional polytope with
 $268$ facets in $8$ orbits, and
  $25448$ vertices in $111$ orbits.  It contains
the exponent vectors of $20992$ monomials in $69$ orbits
which do not appear in $D_{2222}$.
\end{theorem}

In this theorem, the term {\em orbits} refers to the full
symmetry group $B_4$ of the $4$-dimensional cube, which 
has order $384$, and the {\em Newton polytope}
is the convex hull in $\RR^{16}$ of the exponent vectors
of all monomials appearing in $D_{2222}$.

Our computational proof of Theorem~1 is explained in Sections
3 and 4. We list representatives for  each orbit
of facets, vertices, and ``missing monomials,'' and
we discuss key properties found in these data.
Complete data are available at our
supplementary materials website \
\url{bio.math.berkeley.edu/4cube/}.

According to the general theory of Gel'fand, Kapranov and Zelevinsky
\cite[\S 11]{GKZ}, there is a beautiful and deep relationship between the
monomials in the $2 \times 2 \times \cdots \times 2$-hyperdeterminant
and the triangulations of the $n$-dimensional cube. In the language of
polyhedral geometry, this relationship can be stated succinctly as follows:
the {\em secondary polytope} of the $n$-cube equals the Minkowski sum
of the Newton polytopes of the hyperdeterminants corresponding to all the faces of the $n$-cube, up to dimension $n$.
Readers unfamiliar with these geometric concepts will find an elementary
and self-contained introduction in Section 2, where we 
explain the correspondence between the twelve monomials in
$D_{222}$ and the $74$ triangulations of the $3$-cube.

Theorem 1 opened up the possibility of determining the
same correspondence for the $4$-cube. However, the
challenges were formidable, in light of the size of the secondary polytope of the 4-cube:

\begin{theorem}
The $4$-cube has $87959448 $ regular triangulations
in $235277$ orbits. The secondary polytope of the
$4$-cube has $80876$ facets in $334$ orbits.
\end{theorem}

The first statement in Theorem 2 is a known result
due to Pfeifle and Rambau \cite{PR}, who had
computed all regular triangulations of the 4-cube
using Rambau's software {\tt TOPCOM} \cite{Ram}.
We independently verified their enumeration.
The second statement in Theorem 2 is one of the 
contributions  of this paper.

Equipped with Theorems 1 and 2, it took us two more months
of serious computational work to obtain what we consider to be 
the main mathematical result of this paper: the extension of the material
in Section 2 from dimension three to dimension four.
Each regular triangulation of the $4$-cube maps
to a monomial in $D_{2222}$, and the fibers of this
map are the {\em $D$-equivalence classes} of \cite{GKZ}.
We found that, up to symmetry, there are $111$ $D$-equivalence
classes. They are discussed in
Section~6. Section~7 studies the $80876$ facets in Theorem 2, 
and it classifies the coarsest regular polyhedral subdivisions of the $4$-cube.
Most of the algorithms and softwares we developed and used here can easily be generalized and applied to other polynomials and polytopes with symmetry.

Section~5 places our computations  into a larger mathematical
 context. The hyperdeterminant $D_{2222}$ is
the {\em $A$-discriminant} for the matrix of the $4$-cube
$$
A \quad = \quad
\left[
\begin{array}{cccccccccccccccc}
1&1&1&1&1&1&1&1&1&1&1&1&1&1&1&1\\
0&0&0&0&0&0&0&0&1&1&1&1&1&1&1&1\\
0&0&0&0&1&1&1&1&0&0&0&0&1&1&1&1\\
0&0&1&1&0&0&1&1&0&0&1&1&0&0&1&1\\
0&1&0&1&0&1&0&1&0&1&0&1&0&1&0&1
\end{array}
\right].
$$
Dickenstein, Feichtner and Sturmfels \cite{DFS}
recently gave a recipe for computing the
Newton polytope of the $A$-discriminant for any
point configuration $A$. Our computations
are complementary  to their tropical approach.
We represent subdivisions of the cube $A$
by their dual polyhedral complexes, here called
{\em tight spans} as in \cite{Hir}.
The related cases
when $A$ is a second hypersimplex
or a product of two simplices correspond to
finite metrics \cite{SY} and to tropical convexity \cite{DS}
respectively. We also study the irreducible factorizations of all
leading forms of $D_{2222}$ that are supported on
facets of the Newton polytope.

\section{The secondary polytope of the $3$-cube}

As a warm-up for our study of the
$4$-cube, we first discuss our primary
objects of
interest for the $3$-cube. Recall that the {\em Newton polytope} $\mathcal{N}(G)$ of a
polynomial $G$ is the convex hull
of the exponent vectors of the monomials which appear in the
expansion of $G$.  The Newton polytope $\mathcal{N}(D_{222})$
of the hyperdeterminant $D_{222}$ is the convex hull
in $\RR^8$ of the six rows of the matrix
$$ \bordermatrix{
  & x_{000} & x_{001} & x_{010} & x_{011} &   x_{100} & x_{101} & x_{110} &
x_{111} \cr
  c_{000} c_{011} c_{101} c_{110}      & 1 & 0 & 0 & 1 & 0 & 1 & 1 & 0 \cr
   c_{001} c_{010} c_{100} c_{111}   & 0 & 1 & 1 & 0 & 1 & 0 & 0 & 1 \cr
    c_{000}^2c_{111}^2 & 2 & 0 & 0 & 0  & 0 & 0 & 0 & 2 \cr
    c_{001}^2c_{110}^2 & 0 & 2 & 0 & 0 & 0 & 0 & 2 & 0 \cr
    c_{010}^2c_{101}^2 & 0 & 0 & 2 & 0 & 0 & 2 & 0 & 0 \cr
    c_{011}^2c_{100}^2 & 0 & 0 & 0 & 2 & 2 & 0 & 0 & 0 \cr
}.
$$
These six monomials labeling the rows are the {\em extreme monomials} of $D_{222}$,  which means that their exponent vectors are vertices of
the Newton polytope $\mathcal{N}(D_{222})$.
The other six  monomials in $D_{222}$ are not extreme monomials,
since their exponent vectors lie in the relative interiors
of edges of the four-dimensional polytope
  $\mathcal{N}(D_{222})$.
Combinatorially, $\mathcal{N}(D_{222})$
  is a bipyramid over a tetrahedron. Its {\em $f$-vector}
  records the number of faces of dimensions $0$, $1$, $2$, and $3$
respectively:
  \begin{equation}
  \label{smallfvector}
  f\bigl(\mathcal{N}(D_{222})) \quad = \quad (6,14,16,8).
  \end{equation}
The Newton polytope of the $2 \times 2 \times 2$-hyperdeterminant has
  the following irredundant
  presentation by linear equations and inequalities:
\begin{eqnarray*}
& \mathcal{N}(D_{222}) \,\, = \,\,
\{(x_{000},x_{001}, \ldots, x_{111} )\in\mathbb{R}^8 \,:\, \qquad
\qquad \qquad \qquad  \qquad \\
& x_{000}, x_{001}, x_{010}, x_{011}, x_{100}, x_{101}, x_{110}, x_{111}
\geq 0 , \\
& \! x_{000} + x_{001} + x_{010} + x_{011} \,=\,
   x_{000} + x_{001} + x_{100} + x_{101} \,=\, 2 ,\\
& \, \, x_{000} + x_{010} + x_{100} + x_{110} \,=\,
x_{001} + x_{011} + x_{101} + x_{111} \,=\, 2  \,\,
\}.
\end{eqnarray*}
In Section 4, we determine the analogous presentation for $\mathcal{N}(D_{2222})$.

We consider the {\em principal determinant}
of the $3$-cube.
By~\cite[Thm.~10.B.1.2]{GKZ}, this is the
following product of determinants associated to 
faces and vertices:
\begin{eqnarray*}
E_{222} \,\,\, = & D_{222} \cdot
(c_{000} c_{011} - c_{001} c_{010}) \cdot
(c_{000} c_{101} - c_{001} c_{100}) \cdot \\
& (c_{000} c_{110} - c_{010} c_{100}) \cdot
(c_{001} c_{111} - c_{011} c_{101}) \cdot \\
& (c_{010} c_{111} - c_{011} c_{110}) \cdot
(c_{100} c_{111} - c_{101} c_{110}) \cdot \\
& c_{000} \cdot c_{001} \cdot  c_{010} \cdot c_{011} \cdot
c_{100} \cdot c_{101} \cdot c_{110} \cdot c_{111}.
\end{eqnarray*}
The expansion of this polynomial of degree $24$ has
$231$ monomials, of which $74$ are extreme monomials.  The symmetry group
of the $3$-cube, which is the Weyl group $B_3$ of order 48, acts on these $231$ monomials.  The
$74$ extreme monomials come in six orbits:
$$ \begin{matrix}
{\rm  Type} & {\rm Monomial} & \hbox{GKZ Vector} &
\hbox{Orbit Size} \\
1 & -4 c_{000} c_{001}^5 c_{010}^5 c_{011} c_{100}^5 c_{101} c_{110}
c_{111}^5 &
  (1\, 5\, 5\, 1\, 5\, 1\, 1\, 5) & 2 \\
2 &
-c_{000} c_{001}^4 c_{010}^4 c_{011}^3 c_{100}^6 c_{101} c_{110} c_{111}^4 &
  (1 \, 4\, 4\, 3\, 6\, 1\, 1\, 4 ) & 8 \\
3 & c_{000} c_{001}^3 c_{010}^4 c_{011}^4 c_{100}^6 c_{101}^2 c_{110}
c_{111}^3 &
  (1 \, 3 \, 4 \, 4 \, 6 \, 2 \, 1 \, 3) & 24 \\
  4 & c_{000} c_{001}^3 c_{010}^3 c_{011}^5 c_{100}^5 c_{101}^3 c_{110}^3
c_{111} &
  (1 \, 3 \, 3 \, 5 \, 5 \, 3 \, 3 \, 1) & 12 \\
5 & -c_{000} c_{001}^3 c_{010}^3 c_{011}^5 c_{100}^6 c_{101}^2 c_{110}^2
c_{111}^2 &
  (1 \, 3 \, 3 \, 5 \, 6 \, 2 \, 2 \, 2) &  24 \\
6 & c_{000}^2 c_{001}^2 c_{010}^2 c_{011}^6 c_{100}^6 c_{101}^2 c_{110}^2
c_{111}^2 &
  (2 \, 2 \, 2 \, 6 \, 6\, 2 \, 2 \, 2) & 4
\end{matrix}
$$
The Newton polytope $\mathcal{N}(E_{222})$ of the principal
determinant $E_{222}$ is the
{\em secondary polytope} of the $3$-cube.
It is $4$-dimensional and has the  $f$-vector $(74,152,$ $100, 22)$.
A Schlegel diagram of this polytope
was posted by Julian Pfeifle at
\begin{verbatim}
       www.eg-models.de/models/Discrete_Mathematics/
        Polytopes/Secondary_Polytopes/2000.09.031/.
\end{verbatim}

Using {\tt POLYMAKE} \cite{GJ}, we found that
the secondary polytope of the $3$-cube has the following irredundant
  presentation by linear equations and inequalities:
\begin{equation}
\begin{array}{cc}
& \mathcal{N}(E_{222}) \,\, = \,\,
\bigl\{ (x_{000},x_{001}, \ldots, x_{111} )\in\mathbb{R}^8 \,:\, \qquad
\qquad \qquad \qquad  \qquad \\
& x_{000} + x_{001} + x_{010} + x_{011} \,=\,
   x_{000} + x_{001} + x_{100} + x_{101} \,=\, 12 ,\\
&  x_{000} + x_{010} + x_{100} + x_{110} \,=\,
x_{001} + x_{011} + x_{101} + x_{111} \,=\, 12  ,\\
& 1 \, \leq \, x_{000}, \,x_{001}, \,x_{010}, \,x_{011}, \,x_{100},
\,x_{101},
\,x_{110}, \,x_{111}
\,\leq \,6, \\ &
x_{000} + x_{001} \geq 4,\,
x_{000} + x_{010} \geq 4,\,
x_{000} + x_{100} \geq 4 ,\,\\ &
x_{001} + x_{011} \geq 4 ,\,
x_{010} + x_{011} \geq 4,\,
x_{001} + x_{101} \geq 4 \,\bigr\}.
\end{array}
\label{NE222ineq}
\end{equation}

According to the theory of Gel'fand, Kapranov, and Zelevinsky \cite{GKZ},
the $74$ vertices of the secondary polytope $\mathcal{N}(E_{222})$
are in bjiection with the {\em regular triangulations} of the $3$-cube.
For the $3$-cube all triangulations are
regular, so the $3$-cube has $74$ triangulations which come in $6$ orbits.
The vertex of $\mathcal{N}(E_{222})$ corresponding to
a given triangulation $\Pi$ is called the {\em GKZ vector} of $\Pi$.
The $i$th coordinate of the GKZ vector of $\Pi$ is the sum of the
normalized volumes of all tetrahedra in $\Pi$
which contain the $i$th vertex of the cube.
For a general introduction to triangulations, their 
GKZ vectors, and secondary polytopes,
see the forthcoming book by De Loera, Rambau and Santos \cite{DRS}.
The six types of triangulations of the $3$-cube are
depicted in \cite[Figure 1.38]{DRS} and \cite[Figure 3]{DeL}.  Their tight spans, which are polyhedral complexes dual to the triangulations, are shown in Figure 1 below.
Here is a detailed description:

\smallskip

\noindent {\bf Type 1:}  These two triangulations divide the cube into five
tetrahedra. There is
one central tetrahedron of normalized volume two, and four
unimodular tetrahedra which border the central tetrahedron.

\smallskip

The remaining five types of triangulations divide the 3-cube into
six unimodular tetrahedra. Each of these triangulations
uses a diagonal of the cube.

\noindent {\bf Type 2:}
We get these $8$ triangulations  by slicing off the
three vertices adjacent to a fixed vertex.
The remaining bipyramid is cut into three tetrahedra.

\noindent {\bf Type 3:}
These $24$ triangulations also use a diagonal. Of the
other six vertices we pick two that are diagonal on a facet, and we
slice them off.

\noindent {\bf Type 4:}
These twelve triangulations are indexed by ordered pairs of diagonals.
The end points of the first diagonal are sliced off, and the
 remaining octahedron is triangulated using the second diagonal.

\noindent {\bf Type 5:}
These $24$ triangulations are indexed by a diagonal and one other vertex.
That vertex is sliced off, and the remaining polytope is divided into
a pentagonal ring of tetrahedra around the diagonal.

\noindent {\bf Type 6:}
These four triangulations are indexed by the diagonals.
The cube is divided into a hexagonal ring of tetrahedra around the diagonal.

\smallskip

We depict each triangulation of the $3$-cube
by its corresponding tight span, which is
the planar graph dual to the triangulation.
Each vertex of the tight span corresponds to
a tetrahedron in the triangulation, and two
vertices are connected by an edge if the
corresponding tetrahedra share a triangle.  
Regions of the tight span
correspond to interior edges of the triangulation.
The six types of tight spans are shown in Figure~1.  Compare it to Figure 2 in \cite{San}.

\begin{figure}
\includegraphics[scale=0.65]{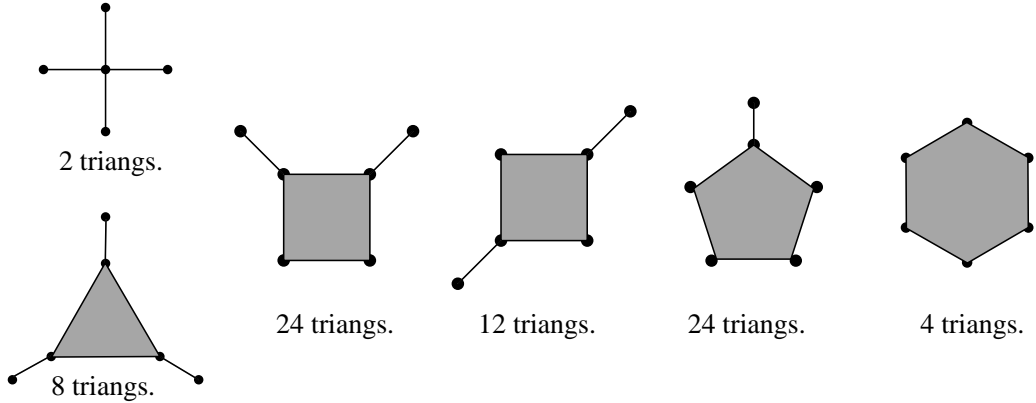}
%\vskip 5cm
\caption{The tight spans dual to the six types of triangulations of the
$3$-cube.
Each tight span represents an orbit of extreme monomials of the principal
determinant $\,E_{222}$.}
\label{cubetriangs}
\end{figure}

The $22$ facets of $\mathcal{N}(E_{222})$ correspond to
proper subdivisions of the $3$-cube which are as coarse as possible.
  The eight facet inequalities
like $x_{000}  \geq 1$ correspond to
slicing off one vertex. The eight inequalities
like $x_{000} \leq 6$ correspond to subdividing
the $3$-cube into three pyramids whose bases
are the square facets disjoint from one vertex.
The six inequalities like $x_{000} + x_{001} \geq 4$
correspond to subdividing the $3$-cube into two triangular prisms.

Since the hyperdeterminant $D_{222}$ is a factor of the principal
determinant $E_{222}$, its Newton polytope
$\mathcal{N}(D_{222})$ is  a {\em Minkowski summand}
of the secondary polytope  $\mathcal{N}(E_{222})$.
This implies the existence of a natural map from the
$74$ vertices of $\mathcal{N}(E_{222})$
onto the six vertices of $\mathcal{N}(D_{222})$.
Equivalently, we have a map from the
regular triangulations of the $3$-cube to
the extreme monomials of $D_{222}$.
A formula for this map can be derived
from Theorem 3.2 in \cite[\S 11.A]{GKZ}.

The extreme monomials of $D_{222}$ come in two $B_3$-orbits.  We now
describe the corresponding
two orbits of {\em $D$-equivalence classes}. Each class consists of all
triangulations which are mapped to a fixed
extreme monomial of $D_{222}$.
The $D$-equivalence class of the monomial $c_{001} c_{010} c_{100} c_{111} $
consists of only one triangulation, namely the triangulation
of type 1 with GKZ vector $(1,5,5,1,5,1,1,5)$.
The $D$-equivalence class of  $c_{000}^2 c_{111}^2$ consists
of all $18$ triangulations which use the diagonal
from $(0,0,0)$ to $(1,1,1)$. Thus the number of triangulations
of the $3$-cube decomposes as follows into the sizes of  the
$D$-equivalence classes:
\begin{equation}
\label{seventyfour} 74 \quad =  \quad 2 \cdot 1 + 4 \cdot 18 \quad =  \quad
1 + 1 + 18 + 18 + 18  + 18.
\end{equation}

We emphasize that all of these results on the 3-cube are
easy and well-known. This section served a purely expository purpose,
namely, to set the stage for the new results on the $4$-cube 
to be presented in the later sections. For instance,
Table \ref{decompose} generalizes the identity (\ref{seventyfour}) 
from the $3$-cube to the $4$-cube.

\section{Schl\"{a}fli's formula and its expansion}

We now consider a multilinear polynomial in four variables $x,y,z,$ and $w$,
\begin{equation*}
\begin{split}
F \,\,= & \,\, c_{0000} + c_{0001} w + c_{0010} z + c_{0011} z w +
c_{0100} y +
c_{0101} y w \\
&+ c_{0110} y z + c_{0111} y z w + c_{1000} x + c_{1001} x w + c_{1010} x
z \\
&+ c_{1011} x z w + c_{1100} x y + c_{1101} x y w + c_{1110} x y z +
c_{1111} x y z w,
\end{split}
\end{equation*}
where the $16$ coefficients $c_{ijkl}$ are
regarded as unknowns. 

\begin{definition}
The {\em hyperdeterminant}
$D_{2222}$ is the unique irreducible polynomial  (up to sign) of content
one in the $16$ unknowns $\,c_{ijkl}\,$ that vanishes whenever the system 
of equations
\begin{equation}
\label{fiveeqns}
F \,\,\, = \,\,\, \frac{\partial F}{\partial x}
\,\,\, = \,\,\, \frac{\partial F}{\partial y}
\,\,\, = \,\,\, \frac{\partial F}{\partial z}
\,\,\, = \,\,\, \frac{\partial F}{\partial w} \quad = \quad 0
\end{equation}
has a solution  $(x_0,y_0,z_0,w_0)$ in $\CC^4$.
\end{definition}

In theory, we can compute $D_{2222}$ by eliminating
the  four variables $x,y,z,w$ from the five equations (\ref{fiveeqns}),
say by using Gr\"obner bases or resultants
\cite{CLO}, but in practice this is infeasible.
The analogous computation for a multilinear polynomial
in three variables, however, is easy to do, and it yields
the expression for $D_{222}$ stated in the Introduction.

Schl\"afli's formula for $D_{2222}$ is obtained as follows.
We replace each of the eight unknowns $\,c_{ijk}\,$
in the $2 \times 2 \times 2$-hyperdeterminant $D_{222}$ by
$\,c_{ijk0} + c_{ijk1} w$. The resulting expression
is a polynomial of degree four in the variable $w$. Its
coefficients are polynomials of degree four in the $c_{ijkl}$. The
discriminant of  this polynomial with respect to the variable $w$ is an expression of degree $24$
in the $c_{ijkl}$. This expression is
Schl\"afli's formula for
the  hyperdeterminant  $D_{2222}$.

  \begin{proposition}
The hyperdeterminant $D_{2222}$
  coincides with the discriminant of
$\,D_{222}( c_{ijk0} + c_{ijk1} w)$,
  considered as a polynomial in $w$,
  divided by $256$.
\end{proposition}

This proposition is a special case of \cite[Theorem 14.4.1]{GKZ}.
Its proof is based on the fact that both polynomials have
degree $24$, which we know for $D_{2222}$ by \cite[Corollary 14.2.10]{GKZ}.
For tensors of larger format (e.g., $2 \! \times \! 2 \! \times\! 2
\!\times \! 2 \!\times \! 2)$
  the same recursive method does not work so well.
  In general, Schl\"afli's formula yields the desired
hyperdeterminant times a large extraneous factor.

The first assertion in Theorem 1 was proved by expanding
Schl\"afli's formula, as follows.  Using ${\tt MAPLE}$ we expressed
$\,D_{222}( c_{ijk0} + c_{ijk1} w)$ as
$$ D_{222}( c_{ijk0} + c_{ijk1} w) \quad = \quad 
b_4 w^4 + b_3 w^3 + b_2 w^2 + b_1 w + b_0 $$
where each $b_i$ is a degree 4 polynomial in the 16 unknowns $c_{ijkl}$.
These expressions for the $b_i$ in terms of the $c_{ijkl}$
were then substituted into the discriminant
\begin{equation}
\label{quarticDisc}
\begin{split}
&
  256 b_0^3 b_4^3
-192 b_0^2 b_1 b_3 b_4^2
-128 b_0^2 b_2^2 b_4^2
+144 b_0^2 b_2 b_3^2 b_4
-27 b_0^2 b_3^4 \\ &
+144 b_0 b_1^2 b_2 b_4^2
-6 b_0 b_1^2 b_3^2 b_4
-80 b_0 b_1 b_2^2 b_3 b_4
+18 b_0 b_1 b_2 b_3^3
+16 b_0 b_2^4 b_4 \\ &
-4 b_0 b_2^3 b_3^2
-27 b_1^4 b_4^2
+18 b_1^3 b_2 b_3 b_4
-4 b_1^3 b_3^3
-4 b_1^2 b_2^3 b_4
+b_1^2 b_2^2 b_3^2.
\end{split}
\end{equation}

After substituting, ${\tt MAPLE}$ was unable to expand and combine the
discriminant's $16$
products due to memory constraints.  Instead, we used ${\tt MAPLE}$ to
substitute and expand each
of the $16$ products (such as $\,b_1^2 b_2^2 b_3^2 \,)$ separately.
We then wrote a
C++ program to merge the $16$ expansions and divide all coefficients by
$256$.
  The result of this merge was a lexicographically sorted list of
  $2894276$ monomials of degree $24$ in the  $16$ unknowns $c_{ijkl}$,
  each with its integer coefficient.  The largest absolute
value of any  coefficient is $ 112464$.

The symmetry group of the $4$-cube is the
Weyl group $B_4$ of order $384$.
It acts on the set of monomials in $D_{2222}$ giving
 $9617$ orbits. 
At our website, the $9617$ lexicographically smallest
monomials in each orbit are listed in lexicographic order.
Each monomial is listed on a separate row, in the format  $\,{\tt
[[ExponentVector], \,Coefficient, \,FaceDimension, \,OrbitSize]}$.
Here ${\tt FaceDimension}$ is the dimension of the smallest face of
$\mathcal{N}(D_{2222})$ containing the exponent vector.
For example, the $192$ monomials in the orbit of
$$ -2 \cdot c_{0011}^2 c_{0101}^2 c_{0110}^3 c_{0111}^5 c_{1000}^7
c_{1001} c_{1010}
c_{1011} c_{1100} c_{1101}
$$
are represented by the row
$\, {\tt [ [0, 0, 0, 2, 0, 2, 3, 5, 7, 1, 1, 1, 1, 1, 0, 0], -2, 3, 192] }$.
The distribution of the $9617$ orbits according to orbit size is as follows:
$$
\begin{array}{cccccccccccc}
2 & 8 & 12 & 16 & 24 & 32 & 48 & 64 & 96 & 192 & 384 \\
\underline{2} & 6 & 6 & 16 & 27 & 24 & 142 & 90 & 577 & 2743 & 5984 \\
\end{array}
$$
The two orbits of size two are represented by the monomial
$$
2008 \cdot c_{0001}^3 c_{0010}^3 c_{0100}^3 c_{0111}^3
c_{1000}^3 c_{1011}^3 c_{1101}^3 c_{1110}^3,
$$
whose exponent vector lies on a $3$-dimensional face, and
the monomial
\begin{eqnarray*} & 112464 \cdot
c_{0000} c_{0001}^2 c_{0010}^2 c_{0011} c_{0100}^2 c_{0101} c_{0110}
c_{0111}^2 \\
& \qquad \quad \, \cdot \,\, c_{1000}^2 c_{1001} c_{1010}
c_{1011}^2 c_{1100} c_{1101}^2 c_{1110}^2 c_{1111},
\end{eqnarray*}
which has the largest coefficient and whose exponent vector lies in the relative interior of $\mathcal{N}(D_{2222})$.  The largest odd coefficient appears in the monomial
$$
-5811 \cdot c_{0001} c_{0010} c_{0011}^4 c_{0100}^2 c_{0101} c_{0110}^2
c_{0111}
c_{1000}^2
c_{1001}^2 c_{1010}^2 c_{1100}^2 c_{1101}^2 c_{1110} c_{1111}
$$
with orbit size $384$ whose exponent vector lies in the relative interior of a $9$-face.

The distribution of the $9617$ orbits according to face dimension is as follows:
$$
\begin{array}{cccccccccccc}
0 & 1 & 2 & 3 & 4 & 5 & 6 & 7 & 8 & 9 & 10 & 11 \\
111 & 230 & 269 & 540 & 1145 & 1862 & 2138 & 1845 & 976 & 405 & 70 & 26
\end{array}
$$

The computation of the above face dimensions and the computation in the rest of this section were all performed after we computed the facet representation of the Newton polytope $\N(D_{2222})$.  In the next section, we present all eight classes of facet inequalities and explain how we obtained them.  

In the remainder of this section we discuss the
$20992$ monomials whose exponent vectors lie in  the Newton polytope
$\mathcal{N}(D_{2222})$ but which have coefficient
zero in the expansion of
hyperdeterminant $D_{2222}$. They come in $69$ $B_4$-orbits.
In the following table we list one representative from each orbit:
\smaller
\begin{equation*}
\begin{array}{cccc}
0 0 0 2 1 1 4 4 4 5 1 0 1 0 1 0 &
0 0 0 2 1 1 4 4 5 4 0 1 1 0 1 0 &
0 0 0 2 1 1 4 4 5 4 1 0 0 1 1 0 &
0 0 0 2 1 1 4 4 6 3 0 1 0 1 1 0 \\
0 0 0 2 1 1 5 3 4 5 1 0 1 0 0 1 &
0 0 0 2 1 1 5 3 5 4 0 1 1 0 0 1 &
0 0 0 2 1 1 5 3 5 4 1 0 0 1 0 1 &
0 0 0 2 1 1 5 3 6 3 0 1 0 1 0 1 \\
0 0 0 2 1 2 2 5 7 1 1 1 1 0 0 1 &
0 0 0 2 1 2 6 1 2 6 2 0 0 1 1 0 &
0 0 0 2 1 3 6 0 2 5 3 0 0 1 0 1 &
0 0 0 2 1 4 4 1 1 4 5 0 1 1 0 0 \\
0 0 0 2 2 1 6 1 2 6 1 1 0 1 1 0 &
0 0 0 2 2 2 3 3 5 3 2 0 0 0 0 2 &
0 0 0 2 2 2 6 0 2 6 2 0 0 0 0 2 &
0 0 0 2 3 1 4 2 3 5 2 0 0 0 0 2 \\
0 0 0 2 4 1 1 4 5 0 1 4 1 1 0 0 &
0 0 0 3 1 1 3 4 4 4 1 0 2 0 1 0 &
0 0 0 3 1 1 4 3 4 4 0 1 2 0 1 0 &
0 0 0 3 1 1 4 3 5 3 0 1 2 0 0 1 \\
0 0 0 3 1 1 4 3 5 3 1 0 0 2 1 0 &
0 0 0 3 1 1 4 3 6 1 1 1 0 3 0 0 &
0 0 0 3 1 1 5 2 6 2 0 1 0 2 0 1 &
0 0 0 3 1 3 4 1 1 4 4 0 2 1 0 0 \\
0 0 0 3 1 4 4 0 5 0 2 2 0 2 0 1 &
0 0 0 3 2 2 2 3 5 0 2 2 1 2 0 0 &
0 0 0 3 2 2 3 2 2 5 2 0 1 0 2 0 &
0 0 0 3 2 2 3 2 4 2 2 1 1 1 0 1 \\
0 0 0 3 2 2 5 0 3 3 1 2 1 1 0 1 &
0 0 0 3 3 2 4 0 4 1 1 3 0 2 0 1 &
0 0 0 3 4 0 1 4 4 1 1 3 2 1 0 0 &
0 0 0 4 1 1 3 3 7 1 0 0 0 2 1 1 \\
0 0 0 4 1 1 4 2 5 2 1 0 0 3 1 0 &
0 0 0 4 1 1 5 1 5 3 0 0 0 2 1 1 &
0 0 0 4 1 1 5 1 6 1 0 1 0 3 0 1 &
0 0 0 4 1 2 5 0 4 3 0 1 1 1 1 1 \\
0 0 0 4 2 1 5 0 4 3 1 0 0 2 0 2 &
0 0 0 5 1 1 4 1 5 1 1 0 0 4 1 0 &
0 0 0 5 1 1 5 0 5 0 1 1 0 5 0 0 &
0 0 0 5 1 1 5 0 6 0 0 1 0 4 0 1 \\
0 0 0 6 1 1 3 1 4 1 1 0 2 3 1 0 &
0 0 0 6 1 1 4 0 6 0 0 0 0 4 1 1 &
0 0 0 6 2 1 2 1 2 2 1 1 4 1 1 0 &
\underline{0 0 1 1 1 1 1 7 7 1 1 1 1 1 0 0} \\
0 0 1 1 1 1 3 5 5 3 1 1 1 1 0 0 &
0 0 1 1 1 2 1 6 7 2 1 0 0 0 1 1 &
0 0 1 1 1 2 1 6 8 1 0 1 0 0 1 1 &
0 0 1 1 1 4 4 1 5 1 1 3 0 1 0 1 \\
0 0 1 1 2 2 2 4 6 2 0 2 0 0 1 1 &
0 0 1 2 1 2 5 1 3 6 0 0 0 0 2 1 &
0 0 1 2 1 4 2 2 4 1 3 1 1 1 0 1 &
0 0 1 2 2 2 0 5 5 1 2 1 2 0 0 1 \\
0 0 1 2 3 1 0 5 4 2 2 1 2 0 0 1 &
0 0 1 2 4 1 0 4 4 1 0 4 2 0 1 0 &
0 0 1 3 1 4 3 0 5 2 1 0 0 0 1 3 &
0 0 1 3 2 2 1 3 3 3 2 0 1 1 2 0 \\
0 0 1 3 2 2 2 2 2 3 2 1 2 1 1 0 &
0 0 1 3 3 3 0 2 5 1 1 1 0 0 2 2 &
0 0 1 4 2 2 0 3 5 1 0 1 1 1 3 0 &
0 0 1 4 2 3 1 1 5 1 0 1 0 1 3 1 \\
0 0 1 4 3 2 0 2 5 2 0 0 0 0 3 2 &
0 0 1 4 3 3 0 1 5 1 1 0 0 0 2 3 &
0 0 1 4 3 3 1 0 5 1 0 1 0 0 2 3 &
0 0 2 3 2 3 1 1 5 2 0 0 0 0 2 3 \\
0 1 1 1 2 2 4 1 3 2 1 3 0 2 1 0 &
0 1 1 2 2 1 4 1 3 2 0 3 1 2 1 0 &
0 1 1 3 2 1 3 1 3 2 2 0 1 2 0 2 &
0 1 1 4 2 1 3 0 3 1 1 1 2 2 0 2 \\
0 1 2 1 2 1 2 3 2 4 2 0 2 0 0 2
\end{array}
\end{equation*}
\normalsize
For instance, the underlined vector
represents the monomial
$$
c_{0010} c_{0011} c_{0100} c_{0101} c_{0110} c_{0111}^7
c_{1000}^7 c_{1001} c_{1010} c_{1011} c_{1100} c_{1101}. $$
The orbit of this monomial has size $96$.
Each of these $69$ monomials listed above does appear during
the expansion of Schl\"afli's formula, namely, when we expand
each of the $16$ terms of the discriminant (\ref{quarticDisc}). But
the coefficients of these monomials sum to zero when we merge to form $D_{2222}$.

The distribution of the $69$ orbits of ``missing monomials''
according to  the face dimension in $\mathcal{N}(D_{2222})$ is as follows:
$$
\begin{array}{cccccccccccc}
0 & 1 & 2 & 3     & 4 & 5   & 6 & 7 & 8 & 9 & 10 & 11 \\
0 & 0 &15 & 3 & 20 & 13 & 7 & 5 & 6 & 0 &  0 &  0 \\
\end{array}
$$

The following method was used to generate all lattice points in
$\mathcal{N}(D_{2222})$ and hence find
the $20992$ ``missing monomials'' in $69$ orbits.  We used the five equations $A \cdot x = (24,12,12,12,12)^T$ and 268 facet defining inequalities of $\N(D_{2222})$ which will be presented in the next section.  By symmetry,
it suffices to generate at least one lattice point in each $B_4$-orbit.  We first listed all
ways to assign eight nonnegative integers summing to $12$ on the first facet, i.e., the first eight coordinates.  There were $50388$ ways in all.  Per the symmetry remark above, we kept  only the $1349$ ways which were lex-min under $B_3$.   We then extended each of these $1349$ facet assignments to the entire $4$-cube as follows.
We fixed two disjoint non-parallel edges on the opposite facet, assigned all possible values for those four entries, and solved for the remaining four entries via the five linear relations above.  Of  the resulting $4$-cube assignments, $87435$ have
non-negative entries.  We then tested these against the facet inequalities and found that $80788$ of them lied in the polytope.  We now had at least one representative from each orbit of lattice points in $\N(D_{2222})$.  From these we found $9686$ orbits  of lattice points. Removing the $9617$ orbits of terms which do appear in $D_{2222}$,
we found precisely the
$ 69$  additional orbits listed above.

\section{Computing the Newton polytope}

In this section we present our census of the vertices, facets, and other faces
of the Newton polytope of the
$2 \times 2 \times 2 \times 2$-hyperdeterminant.
The following result completes the proof of Theorem 1.
Being the analogue of equation (\ref{smallfvector}), it shows the increase
in complexity when passing from the $3$-cube to the $4$-cube.

\begin{proposition}
\label{fD2222}
The $f$-vector of the Newton polytope of $D_{2222}$ equals
\begin{eqnarray*}
  f\bigl(\mathcal{N}(D_{2222})) \quad = &
(\, 25448,\,
178780,\,
555280,\,
1005946,\,
1176976, \\ & \quad
927244,\,
495936,\,
176604, \,39680,\,
5012,\, 268 \, ).
\end{eqnarray*}
\end{proposition}

We computed the vertices of the Newton polytope $\N(D_{2222})$
starting from the exponent vectors of the
2894276 monomials in 9617 orbits described in Section 3.  We took advantage of the symmetry group using the following heuristics:
First, we sorted the set of 9617 lex-min elements, one from each orbit, lexicographically and then chose a small lexicographically contiguous subset.
Within this subset, we removed  points that were redundant
(i.e.,~not vertices)  using {\tt POLYMAKE} \cite{GJ}. 
Whenever a point was found to be redundant, we removed
its entire $B_4$-orbit from our original set of points.
After a few iterations we obtained the convex hull of the 9617 lex-min elements, which has 1794 vertices.  The union of the orbits of those 1794 points contains 484804 points.
We repeated the same process on these 484804 points, finding redundant points in a small subset and removing their whole orbits.  We eventually reached
a subset of $ 25448$ points in $111$ orbits which appeared to be irredundant.
We verified that these points are vertices
by solving $111$ linear programming feasibility problems.  

The {\em tangent cone}  at a vertex $v$ 
is the set $\{v+w \in \RR^{16} : v + \epsilon \cdot w 
\in \N(D_{2222}) \text{ for some } \epsilon > 0\}$.  
The faces of the tangent cone are in 
natural bijection with the faces of 
$\N(D_{2222})$ containing $v$.  
The {\em vertex figure} of a vertex $v$ 
is the $10$-dimensional polytope obtained by cutting off
$v$ from its tangent cone.

\begin{proof}[Proof of Proposition \ref{fD2222}:]
The facets and $f$-vector of $\mathcal{N}(D_{2222})$ were computed as follows:
For each of the $111$ vertex classes, we chose 
a representative $v$ and 
subtracted it from each of the other $25447$ vertices.  
The cone generated by these $25447$ difference vectors is 
(a translate of) the tangent 
cone at $v$.  Using {\tt POLYMAKE} we computed the faces of all dimensions of the $111$ tangent cones.   These computations were fairly fast since $\N(D_{2222})$ is close to being simple.  We now had  at least one representative from the $B_4$-orbit of every face of $\mathcal{N}(D_{2222})$.  Each face is represented by the set of facets
containing that face.  With careful relabeling, we merged the 111 lists of faces and computed their
 orbits, thus obtaining all faces of $\N(D_{2222})$.
\end{proof}

Here is a complete list of representatives for
the $111$ orbits of vertices.  The upper and lower indices attached to each vector
are the coefficient and the orbit size respectively.\label{Dvertices}
\begin{eqnarray*}
  0 0 0 2 0 2 2 6 6 2 2 0 2 0 0 0^{1}_{32} &
  0 0 0 2 0 2 2 6 7 1 1 1 2 0 0 0^{-1}_{192} &
  0 0 0 2 0 2 2 6 8 0 1 1 1 1 0 0^{1}_{192} \\
  0  0  0  2  0  2  2  6  9  0  0  1  0  1  1  0^{-1}_{64} &
  0  0  0  2  0  2  3  5  7  1  0  2  2  0  0  0^{1}_{384} &
  0  0  0  2  0  2  3  5  7  2  1  0  1  0  0  1^{1}_{192}\\
  0  0  0  2  0  2  3  5  8  0  0  2  1  1  0  0^{-1}_{384} &
  0  0  0  2  0  2  3  5  8  1  0  1  1  0  0  1^{-1}_{384} &
  0  0  0  2  0  2  3  5  9  0  0  1  0  1  0  1^{1}_{192} \\
  0  0  0  2  0  2  4  4  8  0  0  2  0  2  0  0^{1}_{192} &
  0  0  0  2  0  2  6  2  2  6  2  0  2  0  0  0^{1}_{96} &
  0  0  0  2  0  2  7  1  2  6  1  1  2  0  0  0^{-1}_{384} \\
  0  0  0  2  0  2  7  1  2  7  1  0  1  0  1  0^{1}_{192} &
  0  0  0  2  0  2  7  1  3  5  0  2  2  0  0  0^{1}_{384} &
  0  0  0  2  0  2  7  1  3  6  0  1  1  0  1  0^{-1}_{384} \\
  0  0  0  2  0  2  7  1  4  5  0  1  0  1  1  0^{1}_{192} &
  0  0  0  2  0  2  8  0  2  6  1  1  1  1  0  0^{1}_{192} &
  0  0  0  2  0  2  8  0  2  7  1  0  1  0  0  1^{-1}_{192} \\
  0  0  0  2  0  2  8  0  3  5  0  2  1  1  0  0^{-1}_{384} &
  0  0  0  2  0  2  8  0  3  6  0  1  1  0  0  1^{1}_{384} &
  0  0  0  2  0  2  8  0  4  5  0  1  0  1  0  1^{-1}_{192} \\
  0  0  0  2  0  3  3  4  7  0  0  3  2  0  0  0^{-1}_{192} &
  0  0  0  2  0  3  3  4  8  1  1  0  0  0  0  2^{1}_{96} &
  0  0  0  2  0  3  3  4  9  0  0  1  0  0  0  2^{-1}_{96} \\
  0  0  0  2  0  3  6  1  3  6  1  0  0  0  2  0^{1}_{384} &
  0  0  0  2  0  3  6  1  4  5  0  1  0  0  2  0^{-1}_{384} &
  0  0  0  2  0  3  7  0  2  5  1  2  2  0  0  0^{1}_{192} \\
  0  0  0  2  0  3  7  0  2  6  1  1  1  0  1  0^{-1}_{384} &
  0  0  0  2  0  3  7  0  3  4  0  3  2  0  0  0^{-1}_{192} &
  0  0  0  2  0  3  7  0  3  5  0  2  1  0  1  0^{1}_{384} \\
  0  0  0  2  0  3  7  0  3  6  1  0  0  0  1  1^{-1}_{384} &
  0  0  0  2  0  3  7  0  4  5  0  1  0  0  1  1^{1}_{384} &
  0  0  0  2  0  3  7  0  4  5  1  0  0  0  0  2^{1}_{192} \\
  0  0  0  2  0  3  7  0  5  4  0  1  0  0  0  2^{-1}_{192} &
  0  0  0  2  0  4  6  0  3  5  1  1  0  0  2  0^{-1}_{384} &
  0  0  0  2  0  4  6  0  4  4  0  2  0  0  2  0^{1}_{384} \\
 \underline{ 0  0  0  2  0  4  6  0  6  0  0  4  0  2  0  0}^{16}_{192} &
 \underline{  0  0  0  2  0  5  5  0  5  0  0  5  2  0  0  0}^{-16}_{96} &
  \underline{ 0  0  0  2  0  5  5  0  7  0  0  3  0  0  0  2}^{-16}_{96} \\
  0  0  0  2  1  1  1  7  7  1  1  1  2  0  0  0^{1}_{96} &
  0  0  0  2  1  1  1  7  7  1  2  0  1  1  0  0^{1}_{192} &
  0  0  0  2  1  1  1  7  8  0  1  1  1  1  0  0^{-1}_{384} \\
  0  0  0  2  1  1  1  7  8  1  1  0  0  1  1  0^{-1}_{192} &
  0  0  0  2  1  1  1  7  9  0  0  1  0  1  1  0^{1}_{192} &
  0  0  0  2  1  1  2  6  7  1  0  2  2  0  0  0^{-1}_{384} \\
  0  0  0  2  1  1  2  6  7  1  2  0  0  2  0  0^{-1}_{384} &
  0  0  0  2  1  1  2  6  7  2  1  0  1  0  0  1^{-1}_{384} &
  0  0  0  2  1  1  2  6  8  0  0  2  1  1  0  0^{1}_{384} \\
      0  0  0  2  1  1  2  6  8  0  1  1  0  2  0  0^{1}_{384} &
  0  0  0  2  1  1  2  6  8  1  0  1  1  0  0  1^{1}_{192} &
  0  0  0  2  1  1  2  6  8  1  1  0  0  1  0  1^{1}_{384} \\
   0  0  0  2  1  1  2  6  9  0  0  1  0  1  0  1^{-1}_{192} &
  0  0  0  2  1  1  7  1  1  7  1  1  2  0  0  0^{1}_{96} &
  0  0  0  2  1  1  7  1  1  7  2  0  1  1  0  0^{1}_{192} \\
  0  0  0  2  1  1  7  1  1  8  1  0  1  0  1  0^{-1}_{384} &
  0  0  0  2  1  1  7  1  2  6  0  2  2  0  0  0^{-1}_{192} &
  0  0  0  2  1  1  7  1  2  6  2  0  0  2  0  0^{-1}_{384} \\
    0  0  0  2  1  1  7  1  2  7  0  1  1  0  1  0^{1}_{384} &
  0  0  0  2  1  1  7  1  2  7  1  0  0  1  1  0^{1}_{384} &
  0  0  0  2  1  1  7  1  3  6  0  1  0  1  1  0^{-1}_{384} \\
    0  0  0  2  1  1  8  0  1  8  1  0  1  0  0  1^{1}_{192} &
  0  0  0  2  1  1  8  0  2  7  0  1  1  0  0  1^{-1}_{384} &
  0  0  0  2  1  1  8  0  2  7  1  0  0  1  0  1^{-1}_{384} \\
  0  0  0  2  1  1  8  0  3  6  0  1  0  1  0  1^{1}_{384} &
  0  0  0  2  1  2  2  5  7  0  0  3  2  0  0  0^{1}_{192} &
  0  0  0  2  1  2  2  5  8  1  1  0  0  0  0  2^{-1}_{96}\\
    0  0  0  2  1  2  2  5  9  0  0  1  0  0  0  2^{1}_{48} &
  0  0  0  2  2  0  1  7  6  3  1  0  1  0  1  0^{-1}_{384} &
  0  0  0  2  2  0  1  7  7  2  0  1  1  0  1  0^{1}_{192} \\
  0  0  0  2  2  0  1  7  7  2  1  0  0  1  1  0^{1}_{384} &
  0  0  0  2  2  0  2  6  6  2  0  2  2  0  0  0^{1}_{96} &
  0  0  0  2  2  0  2  6  6  2  2  0  0  2  0  0^{1}_{192} \\
                                              \end{eqnarray*} \begin{eqnarray*}
  0  0  0  2  2  0  2  6  6  3  1  0  1  0  0  1^{1}_{384} &
  0  0  0  2  2  0  2  6  7  2  1  0  0  1  0  1^{-1}_{384} &
  0  0  0  2  2  1  7  0  2  7  1  0  0  0  0  2^{1}_{96} \\
   0  0  0  2  2  1  7  0  3  6  0  1  0  0  0  2^{-1}_{192} &
  0  0  0  2  3  0  0  7  5  2  3  0  1  1  0  0^{1}_{384} &
  0  0  0  2  3  0  0  7  6  1  2  1  1  1  0  0^{-1}_{384} \\
  0  0  0  2  3  0  0  7  6  2  2  0  0  1  1  0^{-1}_{192} &
  0  0  0  2  3  0  0  7  7  1  1  1  0  1  1  0^{1}_{192} &
  0  0  0  2  3  0  1  6  6  3  1  0  0  0  1  1^{1}_{384} \\
  0  0  0  2  3  0  2  5  6  3  1  0  0  0  0  2^{-1}_{192} &
  0  0  0  2  3  0  6  1  3  6  0  1  0  0  0  2^{1}_{96} &
  0  0  0  2  4  0  0  6  5  2  3  0  0  1  0  1^{-1}_{384} \\
  0  0  0  2  4  0  0  6  6  1  2  1  0  1  0  1^{1}_{192} &
  0  0  0  2  4  3  3  0  5  0  0  5  0  0  0  2^{-1}_{96} &
  0  0  0  2  5  0  0  5  5  2  2  1  0  0  0  2^{1}_{48} \\
  0  0  0  3  0  3  3  3  9  0  0  0  0  0  0  3^{1}_{16} &
  0  0  0  3  0  3  5  1  4  5  0  0  0  0  3  0^{1}_{192} &
  0  0  0  3  0  3  6  0  2  6  1  0  1  0  2  0^{1}_{96} \\
  0  0  0  3  0  3  6  0  4  5  0  0  0  0  2  1^{-1}_{192} &
  0  0  0  3  0  3  6  0  6  3  0  0  0  0  0  3^{-1}_{64} &
  0  0  0  3  0  4  5  0  3  5  1  0  0  0  3  0^{1}_{192} \\
  0  0  0  3  1  1  1  6  7  1  1  0  2  0  0  1^{-1}_{192} &
  0  0  0  3  1  1  1  6  8  0  1  0  1  1  0  1^{1}_{192} &
  0  0  0  3  1  1  1  6  9  0  0  0  0  1  1  1^{-1}_{64} \\
  0  0  0  3  1  1  7  0  1  7  1  0  1  1  1  0^{1}_{96} &
  0  0  0  3  1  1  7  0  2  6  1  0  0  2  1  0^{-1}_{192} &
  0  0  0  3  1  1  7  0  3  6  0  0  0  1  1  1^{-1}_{192} \\
  0  0  0  4  0  4  4  0  4  0  0  4  4  0  0  0^{1}_{96} &
  \underline{0  0  0  4  0  4  4  0  8  0  0  0  0  0  0  4}^{-27}_{16} &
  0  0  0  8  1  1  1  1  1  1  1  1  8  0  0  0^{1}_{48} \\
  0  0  0  9  1  1  1  0  1  1  1  0  7  1  1  0^{1}_{48} &
  0  0  0  9  1  1  1  0  1  1  1  0  8  0  0  1^{-1}_{96} &
  0  0  1  1  1  1  0  8  8  1  0  1  1  0  1  0^{1}_{92} \\
  0  0  1  1  1  1  1  7  7  2  1  0  1  0  0  1^{1}_{192} &
  0  0  1  1  1  1  1  7  8  1  0  1  1  0  0  1^{-1}_{192} &
  0  0  1  1  1  2  0  7  7  1  2  0  1  0  0  1^{1}_{192} \\
  0  0  1  1  1  2  0  7  8  1  1  0  0  0  1  1^{-1}_{192} &
\underline{0  1  1  0  1  0  0  9  9  0  0  1  0  1  1  0}^{1}_{8} &
  0  1  1  0  1  0  6  3  3  6  0  1  0  1  1  0^{1}_{32}
\end{eqnarray*}

The distribution of the $111$ orbit sizes is
$$
\begin{array}{cccccccc}
8 & 16 & 32 & 48 & 64 & 96 & 192 & 384 \\
1 & 2   & 2   & 4   & 3 & 16 & 44 & 39 \\
\end{array}
$$
The distribution of the $111$ corresponding extreme monomial
coefficients is
$$
\begin{array}{ccccc}
1 & -1 & 16 & -16 & -27 \\
60 & 47 & 1 & 2 & 1 \\
\end{array}
$$
All but four of the monomials have coefficient $\pm 1$.
The largest absolute value of a
coefficient appears in  the monomial
$-27 \cdot c_{0011}^4 c_{0101}^4 c_{0110}^4 c_{1000}^8 c_{1111}^4$ 
which has orbit size 16.  Its  $D$-equivalence class will 
be discussed in Example \ref{fatclass}.  The monomial
$16 \cdot c_{0011}^2 c_{0101}^4 c_{0110}^6 c_{1000}^6 c_{1011}^4 c_{1101}^2$
has orbit size $192$, and there are two monomials with
coefficient $-16$, each having orbit size $96$.  Their $D$-equivalence classes will be discussed in Example \ref{coeff16}.

Among the $111$ types of vertices, $35$ are simple.
The last underlined vertex, representing the monomial $8 \cdot c_{0001} c_{0010} c_{0100} c_{0111}^9 c_{1000}^9 c_{1011} c_{1101} c_{1110}$, has the largest vertex figure.
That vertex figure has the  $f$-vector
$$(67,873,4405,11451,17440,16452,9699,3446,667,56).$$
Thus this distinguished vertex is adjacent to
$67$ other vertices, and it lies in $56$ facets.
In Example \ref{largestclass} we shall see that this vertex
corresponds to the largest $D$-equivalence class
of triangulations.
Table 1 shows the distribution of the $111$ types of vertices
  according to the number of incident edges and facets.

\begin{table}
\footnotesize
\begin{tabular}{|c|c|c|c|c|c|c|c|c|c|c|c|c|c|c|c|c|c|c|c|c|c|c|}
\hline
       &11&12&13&14&15&16&17&18&19&20&21&22&26&28&29&39&56&\!total\!\\ \hline
  11&35&     &    &     &     &     &    &     &     &    &     &     &
&     &
  &     &    & 35\\ \hline
12 &     &14&    &     &     &     &    &     &     &    &     &     &
&     &
  &     &    & 14\\ \hline
13 &     &13& 5 & 1 &     &     &    &     &     &    &     &     &    &
   &     &
     &    & 19\\ \hline
14 &     &     & 1 & 2 &    &     &    &     &     &    &     &     &    &
     &
&     &    & 3\\ \hline
15 &     &1 & 3 & 1 & 1 &      & 1 &    &     &    &     &     &    &
&     &
  &    & 7 \\ \hline
16 &     &     & 1 & 1 & 2 &     &    &     &     &    &     &     &    &
    &
&     &    &  4 \\ \hline
17 &     &     & 1 & 1 & 3 & 1 & 1  &     &     &    &     &     &    &
  &     &
    &    &  7\\ \hline
18 &     &     &    &     & 1 & 3 &     & 1 & 1  &    &     &     &    &
   &     &
     &    & 6 \\ \hline
19 &     &     &    & 1  &     &     &1 &     &     &    &     &     &
&     &
  &     &    & 2 \\ \hline
20 &     &     &    &     &     &     & 1 &     &     &    &     &     &
  &     &
    &     &    & 1 \\ \hline
21 &     &     &    &     &     &     &    &     & 1  &1  &1   &     &
&     &
  &     &    &3  \\ \hline
22 &     &     &    &     &     &     &    &     &     &    &     &   1 &
   &     &
     &     &    & 1 \\ \hline
23 &     &     &    &     & 1  &     &    &     &     &    &     &     &
  &     &
    &     &    & 1 \\ \hline
25 &     &     &    &     &     &     &    &     &     &    & 1  &     & 1
&     &
   &     &    &2  \\ \hline
26 &     &     &    &     &     &     &    &     &     &    &     &     &
   &   1&
    &     &    & 1 \\ \hline
27 &     &     &    &     &     &     &    &     &     &    &     & 1  &
  &     &
    &     &    & 1 \\ \hline
30 &     &     &    &     &     &     &    &     &     &    &     &     &
   &   1 &
     &     &    & 1 \\ \hline
32 &     &     &    &     &     &     &    &     &     &    &     &     &
   &     &
  1  &     &    & 1 \\ \hline
42 &     &     &    &     &     &     &    &     &     &    &     &     &
   &     &
     &  1 &    & 1 \\ \hline
67 &      &     &    &     &     &     &    &     &     &    &     &     &
    &
&     &    & 1 & 1 \\ \hline
\!total\!&35&28&11&7&8&4&4&1&2&1&2&2&1&2&1&1&1 & 111 \\ \hline
\end{tabular}
\vskip .2cm
\caption{Distribution of the 111 vertices
according to the number of facets (columns) and
edges (rows) containing the vertex.}
\label{tab:vertices}
\end{table}

The $268$ facets of  $\N(D_{2222})$
come in eight orbits. In Table 2
  we present the $f$-vector
and the orbit size for each facet.
We discuss the geometry of
each of the eight facet types.
Since $\N(D_{2222})$ is a Minkowski summand of the
secondary polytope $\N(E_{2222})$, each
facet normal of $\N(D_{2222})$ is a facet normal
of $\N(E_{2222})$ and hence corresponds to
a {\em coarsest regular subdivision} of the $4$-cube.
We depict each subdivision by its {\em tight span}, which
is the complex of bounded faces of the polyhedron
$\, \{ u \in \RR^5 : u \cdot A \leq w \}$,
where  $A$ is the $5 \times 16$-matrix
in Section 1 and $w \in \RR^{16}$ is the facet normal in question.

\begin{table}
\begin{tabular}{|c|c|c|}
\hline
  & \text{$f$-vector} & \!\text{orbit}\! \\ \hline
1& \!(11625,72614,197704,308238,303068,194347,80874,20906,3021,187)\!\!&
 16 \\ \hline
2&  (4112, 25746, 71456,115356,119228,81590,36802,10488, 1704, 122)& 12 \\ \hline
3& (363, 2289, 6538, 10996, 11921, 8581, 4080, 1239, 225, 22) & 8  \\ \hline
4& (938, 5226, 13182, 19716, 19263, 12765, 5758, 1721, 318, 31) & 32  \\ \hline
5& (336, 1937, 5126, 8121, 8468, 6022, 2928, 950, 194, 22)& 48 \\ \hline
6& (289, 1624, 4228, 6636, 6894, 4914, 2413, 798, 168, 20)   & 96 \\ \hline
7& (450, 2526, 6522, 10103, 10315, 7195, 3440, 1099, 220, 24)  & 48\\ \hline
8& (681, 3906, 10323, 16407, 17194, 12264, 5933, 1877, 357, 34) & 8 \\ \hline
\end{tabular}
\vskip .2cm
\label{tab:facets}
\caption{The eight types of facets of the Newton polytope $\N(D_{2222})$}
\end{table}

\medskip

\noindent {\bf Facet 1:}
These $16$ facets are defined by inequalities like
$\,x_{0000} \geq 0 $. The tight span of the corresponding
subdivision is a line segment.
The two maximal cells are the simplex obtained by slicing off one vertex
and the convex hull of the other 15 vertices, which has $f$-vector $(15,34,28,9)$.

\smallskip
\noindent {\bf Facet 2:}
These $12$ facets are given by inequalities like
$$
x_{0000}+x_{0001}+x_{0010}+x_{0011} \geq 2.
$$
  The tight span is a line segment. The two maximal cells
in the corresponding subdivision
  have $f$-vector $(12,24,19,7)$.  Here we
cut the 4-cube with the hyperplane containing the eight vertices
in a pair of opposite 2-faces.

\smallskip
\noindent {\bf Facet 3:}
These $8$ facets are given by inequalities like
$$
2 \cdot x_{0000}+x_{0001}+x_{0010}+x_{0100}+x_{1000} \geq 2.
$$
  The tight span is a line segment, and the $f$-vector of both
maximal cells is $(11,28,26,9)$.
Here we cut the 4-cube with the hyperplane spanned
by the six vertices which have coordinate sum two, whose convex hull is an octahedron.

\smallskip
\noindent {\bf Facet 4:}
These $32$ facets are given by inequalities like
$\, x_{0000}+x_{0001} \leq 9 $.
The tight span is a triangle.  Each of the three maximal cells is a prism
over a square pyramid and has $f$-vector $(10,21,18,7)$.  They are
formed by fixing an edge of the $4$-cube and joining that
edge to the three facets which are disjoint from that edge.
Any two cells intersect in a triangular prism, and all
three cells intersect in the rectangle formed by the fixed
edge and its opposite edge.

\begin{figure}
\includegraphics[scale=0.5]{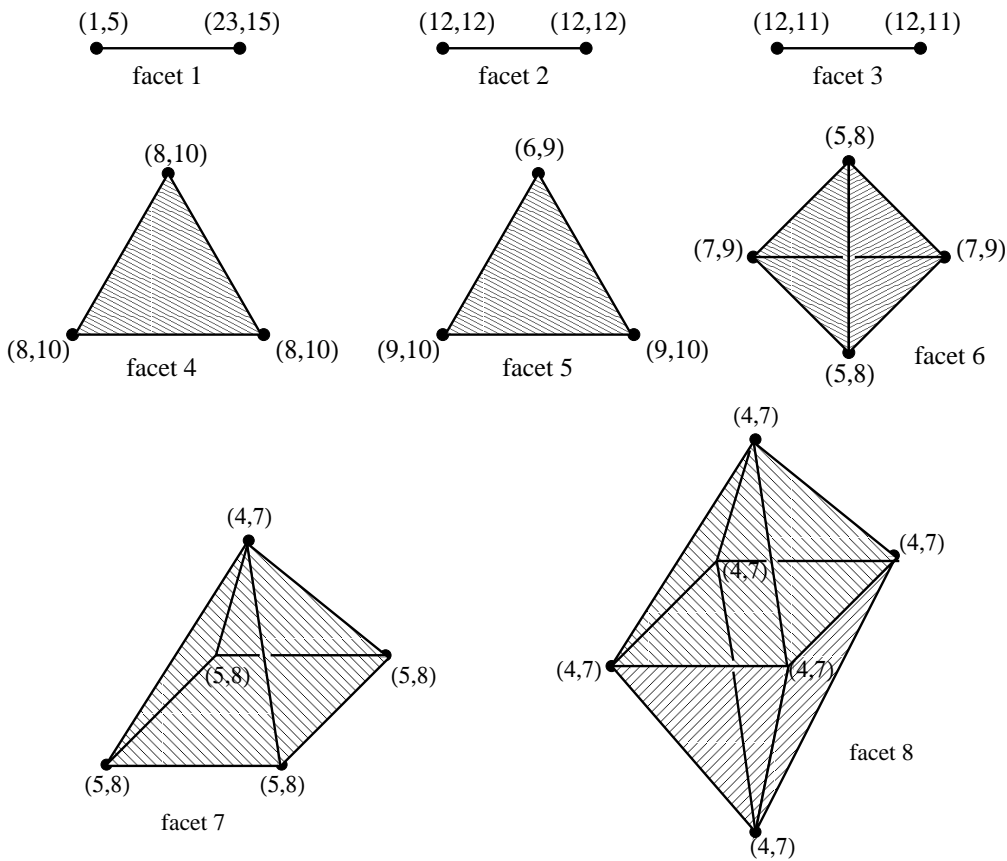} 
%\vskip 6cm
\caption{The tight spans dual to the subdivisions
of the 4-cube  corresponding to the
facets of $\N(D_{2222})$.  Each vertex  of a tight span
is labeled by the normalized volume and the number of vertices of the
corresponding
maximal cell in the subdivision.}
\end{figure}

\smallskip
\noindent {\bf Facet 5:}
These 48 facets are given by inequalities like
$$
x_{0000}+x_{0001}+x_{0010}+x_{0100}+x_{0110}+x_{1000}+x_{1001} \geq 3 .
$$
The tight span is a triangle.  Two of the three maximal cells have $f$-vector
$(10,23,21,8)$  and intersect in a square pyramid; the other cell has
$f$-vector
$(9,18,15,6)$ and intersects each of the first two in a triangular prism.
The intersection of all three cells is a square.

\smallskip
\noindent {\bf Facet 6:}
These 96 facets are given by inequalities like
$$
x_{0000}+x_{0001}+x_{0010}+x_{0011}+x_{0100}+x_{0110}+x_{1000}+x_{1001}
\geq 4.
$$
The tight span is a tetrahedron.  Two of the maximal cells have
$f$-vector $(8,18,17,7)$ and the other two have $f$-vector
$(9,21,20,8)$. The intersection of all four cells is a segment which is a
diagonal of the $4$-cube.

\smallskip
\noindent {\bf Facet 7:}
These 48 facets are given by inequalities like
$$
x_{0000}+x_{0001}+x_{0010}+x_{0011}+x_{0100}+x_{1000}+x_{1100} \leq 19.
$$
The tight span is a square pyramid.  
The maximal cell corresponding to the tip of the pyramid is the convex hull of two squares in complementary dimensions that share a vertex.  It has $f$-vector $(7,17,18,8)$ and meets each of the other four cells in a triangle. 
The other four maximal cells have $f$-vector $(8,18,17,7)$, and their intersection is a triangle.   All the cells meet along a diagonal segment.

\smallskip
\noindent {\bf Facet 8:}
These 8 facets are given by inequalities like
$$
2 \cdot x_{0000}+x_{0001}+x_{0010}+x_{0100}+ x_{1000} \leq 18.
$$
These facets are indexed by diagonals.  We can identify
the $4$-cube with a boolean lattice when we fix a diagonal.  The maximal
antichain in this lattice has the structure of an octahedron.  Each point in the
antichain uniquely determines a maximal cell in the subdivision.  Hence the tight span is also an octahedron.  Each maximal cell is the convex hull of two squares in complementary dimensions joined
at a vertex, and each has $f$-vector $(7,17,18,8)$.

\smallskip

This concludes our discussion of the facets of
the Newton  polytope $\mathcal{N}(D_{2222})$.
The $B_4$-orbits of the eight
stated inequalities yield $268$ facet inequalities. These together with the identity
$\, A \cdot x \, = \, (24,12,12,12,12)^T \,$
gives an irredundant  presentation
of  $\mathcal{N}(D_{2222})$
  by linear equations and inequalities.

\medskip

\section{Tight spans and  $A$-discriminants of other
$0/1$-polytopes}

The concepts and computations presented in this paper
make sense for any $d \times n$-integer matrix $A$ whose row span
contains the all-ones vector $(1,1,\ldots,1) $. Following
\cite{GKZ}, the $A$-discriminant $D_A$ is an irreducible factor of 
the principal $A$-determinant $E_A$.  The
vertices of the secondary polytope $\mathcal{N}(E_A)$
correspond to regular triangulations of $A$, and these
map onto the vertices of the Newton polytope
$\mathcal{N}(D_A)$. The fibers of
that map (i.e., the $D$-equivalence classes) and the 
normal fan of $\mathcal{N}(D_A)$ have
recently been characterized in \cite{DFS}.

For fixed $A$ and an arbitrary row vector $w \in \RR^n$,
we consider the polyhedron
$$ P_w \quad = \quad \{\,u \in \RR^d \,: \, u \cdot A \leq w \,\}. $$
This polyhedron is dual to the regular triangulation of $A$
defined by $w$. The {\em tight span} of $(A,w)$ is the
complex of bounded faces of the polyhedron $P_w$. For
experts in tropical geometry, we note that the 
$(d-2)$-skeleton of $P_w$ is the {\em tropical hypersurface} of $(A,w)$.
The polyhedron $P_w$ also occurs naturally when
relating regular triangulations to linear programming 
duality as in \cite[\S 1.2]{DRS}.

We propose that, for many matrices $A$ of
interest in combinatorics and discrete convexity  \cite{Hir},
the tight span is an excellent geometric 
representation of the data $(A,w)$.
In Section 2 we studied the case when $A$ is the
$4 \times 8$-matrix representing 
the $3$-cube, and in Sections 4, 6 and 7
we are concerned with the case when $A$ is the
$5 \times 16$-matrix representing the $4$-cube. 
In what follows we present two
related situations which have appeared
in the recent literature.

\begin{example}[\bf Tropical Polytopes]
\label{ex1} \rm
Let $A$ be the $(r+s) \times (r \cdot s)$-matrix
which represents the direct product of an
$(r-1)$-simplex $\Delta_{r-1}$
with an $(s-1)$-simplex $\Delta_{s-1}$.
The principal $A$-determinant $E_A$
is the product of all subdeterminants
(of all sizes) of an $r \times s$-matrix 
of unknowns $(c_{ij})$; see \cite[page 303]{GKZ}.

For instance, if $r=s=3$ then $A$ represents the product of two triangles,
\begin{equation}
\label{triangletimestriangle}
A \quad = \quad
\left[ 
\begin{array}{ccccccccc}
1 & 1 & 1 & 0 & 0 & 0 & 0 & 0 & 0 \\
0 & 0 & 0 & 1 & 1 & 1 & 0 & 0 & 0 \\
0 & 0 & 0 & 0 & 0 & 0 & 1 & 1 & 1 \\
1 & 0 & 0 & 1 & 0 & 0 & 1 & 0 & 0 \\
0 & 1 & 0 & 0 & 1 & 0 & 0 & 1 & 0  \\
0 & 0 & 1 & 0 & 0 & 1 & 0 & 0 & 1 
\end{array}
\right],
\end{equation}
and the $A$-discriminant $D_A$ is the determinant of the matrix
$$  \begin{pmatrix} c_{11} & c_{12} & c_{13}\\
c_{21} & c_{22} & c_{23}\\
c_{31} & c_{32} & c_{33}
\end{pmatrix}.$$
The degree $30$ polynomial $E_A$ is the product
of all $19$ subdeterminants of this matrix. 
The extreme monomials of $E_A$ appear in
\cite[Figure 39, page 250]{GKZ}.

The Newton polytope $\mathcal{N}(E_A)$ of the product of all
subdeterminants is the secondary polytope whose faces
correspond to regular polyhedral subdivisions of
$\Delta_{r-1} \times \Delta_{s-1}$. The corresponding
tight spans are precisely the {\em tropical polytopes},
which are obtained as tropical convex hulls \cite{DS}
of $r$ points in $(s-1)$-space. The combinatorial classification of
such polytopes is the analogue of what will be
accomplished for the $4$ cube in
the next section. We refer to recent work of Santos \cite{San}
for a discussion. His
Section 5 deals with tropical polytopes, 
and his Figure 2 is essentially the same as our Figure 1.
Our number  $235277$ in Theorem 2 is the $4$-cube analogue 
of the numbers in the table at the end of
Section 4 in \cite{DS} including the number $35$ in  
\cite[Fig.~1]{DeL} and in \cite[Fig.~6]{DS}.
\qed
\end{example}

\begin{example} [\bf Injective hulls of finite metric spaces]
 \label{ex2}
The term ``tight span'' is derived from the special case when
 $A$ represents the {\em second hypersimplex}.
As shown in \cite{Hir, SY}, here the
tight spans are precisely the {\em injective hulls
of finite metric spaces}, which play an important
role in phylogenetic combinatorics. 
The case studied in \cite{SY}
concerns metrics on six points, where
$$ A \quad = \quad
\left[ 
\begin{array}{ccccccccccccccc}
1 & 1 & 1 & 1 & 1 & 0 & 0 & 0 & 0 & 0 & 0 & 0 & 0 & 0 & 0 \\
1 & 0 & 0 & 0 & 0 & 1 & 1 & 1 & 1 & 0 & 0 & 0 & 0 & 0 & 0 \\
0 & 1 & 0 & 0 & 0 & 1 & 0 & 0 & 0 & 1 & 1 & 1 & 0 & 0 & 0 \\
0 & 0 & 1 & 0 & 0 & 0 & 1 & 0 & 0 & 1 & 0 & 0 & 1 & 1 & 0 \\
0 & 0 & 0 & 1 & 0 & 0 & 0 & 1 & 0 & 0 & 1 & 0 & 1 & 0 & 1 \\
0 & 0 & 0 & 0 & 1 & 0 & 0 & 0 & 1 & 0 & 0 & 1 & 0 & 1 & 1 \\
\end{array}
\right] .
$$
Here the $A$-discriminant $D_A$ is the determinant of the
symmetric matrix
$$  \begin{pmatrix} 
0 & c_{12} & c_{13} & c_{14} & c_{15} & c_{16} \\
c_{12} & 0 & c_{23} & c_{24} & c_{25} & c_{26} \\
c_{13} & c_{23} & 0 & c_{34} & c_{35} & c_{36} \\
c_{14} & c_{24} & c_{34} & 0 & c_{45} & c_{46} \\
c_{15} & c_{25} & c_{35} & c_{45} & 0 & c_{56} \\
c_{16} & c_{26} & c_{36} & c_{46} & c_{56} & 0 \\
\end{pmatrix}.
$$
Up to a constant, the principal $A$-determinant $E_A$ is
the product of all principal
minors of size $\geq 3$. This is a polynomial
of degree $156$ having $194160$ extreme terms in 
$339$ symmetry classes \cite[Theorem 1]{SY}.
Thus our number  $235277$ is the $4$-cube analogue 
to the number $339$ of generic  six-point metrics.
\qed
\end{example}

The prominent role of the second hypersimplex
in phylogenetic combinatorics
raises the question of what happens for other
hypersimplices. We propose the following 
 problem for further
mathematical and computational research.

\begin{problem} [\bf The hypersimplex $\Delta(6,3)$]
\rm
 \label{ex3}  \hfill \break
Let $A$ be the $6 \times 20$-matrix whose
columns are the vertices
$e_i+e_j+e_k$ ($1 \leq i < j < k \leq 6$) of
the third hypersimplex $\Delta(6,3)$.
What is the degree of the $A$-discriminant $D_A$~?
Can the  $14$-dimensional polytope $\N(D_A)$ be computed~?
Can the monomial expansion of the $A$-discriminant $D_A$ 
be computed~?
\end{problem}

All matrices $A$ in the examples above have their entries
in $\{0,1\}$, so they represent subpolytopes
of cubes of appropriate dimensions.
In what follows we examine the $A$-discriminants of
various subpolytopes of the $4$-cube. 
These arise naturally as irreducible factors in the initial forms
of the hyperdeterminant.

For any vector $w$ in $\RR^{16}$ we define the
{\em initial form} $\,{\rm in}_w(D_{2222})\,$
as the sum of all terms in $D_{2222}$ having  minimal $w$-weight.
If $w$ is generic then  $\,{\rm in}_w(D_{2222})\,$
is a monomial, and we have classified all of these
monomials in Section 4. We now consider the other
extreme case when $w$ is as non-generic as possible.
More precisely, we pick $w$ among the normal vectors 
to the facets
of $\mathcal{N}(D_{2222})$.

\begin{proposition}  [\bf All maximal initial forms of the hyperdeterminant]
 \label{all max}
The following list is the classification of the initial forms
of hyperdeterminant $D_{2222}$ corresponding to all facets of its
Newton polytope. The eight symmetry classes of facets
of $\N(D_{2222})$ are listed
in the same order as in Section~4. \rm

\smallskip

\noindent {\bf Facet 1:}
% $$e_{0000} \geq 0.$$ 
If $w = e_{0000} $ then
${\rm in}_w(D_{2222})$ is obtained from
$D_{2222}$ by setting $c_{0000}=0$.
This initial form is irreducible: It is the
$A$-discriminant where $A$ is the
$15$-point configuration  
obtained from the $4$-cube by removing 
one vertex.

\smallskip
\noindent {\bf Facet 2:}
% $$ e_{0000}+e_{0001}+e_{0010}+e_{0011} \geq 2. $$
If $w = e_{0000}+e_{0001}+e_{0010}+e_{0011}$ then
${\rm in}_w(D_{2222})$ has $67230$ terms and 
factors as
$c_{1010} c_{1001} - c_{1000} c_{1011}$ times
$c_{0101} c_{0110} - c_{0100} c_{0111}$ times
the square of the $2 \times 2 \times 2$-hyperdeterminant of 
$[c_{1000}, c_{1001}, c_{1010}, c_{1011}, c_{0100} ,c_{0101} , c_{0110} ,c_{0111}]
$
times the product of two larger factors,
each having $66$ terms of degree six,
which are the $A$-discriminants of 
the two maximal cells in the subdivision.

\smallskip
\noindent {\bf Facet 3:}
% $$  2 \cdot e_{0000}+e_{0001}+e_{0010}+e_{0100}+e_{1000} \geq 2.$$
If $ w =  2 \cdot e_{0000}+e_{0001}+e_{0010}+e_{0100}+e_{1000}$ then
${\rm in}_w(D_{2222})$ equals
\begin{smaller}
\begin{eqnarray*}
& (1/4) \cdot
c_{0011}
c_{0101}
c_{0110}
c_{1001}
c_{1010}
c_{1100} \cdot
\biggl\{
{\rm det} \left[ \begin{matrix}
     0      & c_{0011} & c_{0101} & c_{1001} \\
   c_{0011} &    0     & c_{0110} & c_{1010} \\
   c_{0101} & c_{0110} &     0    & c_{1100} \\
   c_{1001} & c_{1010} & c_{1100} &   0   
   \end{matrix} \right] \biggr\}^{\mathbf 2} \times \\ &
\!\!\!
{\rm det} \! \left[ \begin{matrix}
 2 c_{0000} & c_{0001} & c_{0010} & c_{0100} & c_{1000} \\
   c_{0001} &   0      & c_{0011} & c_{0101} & c_{1001} \\
   c_{0010} & c_{0011} &    0     & c_{0110} & c_{1010} \\
   c_{0100} & c_{0101} & c_{0110} &     0    & c_{1100} \\
   c_{1000} & c_{1001} & c_{1010} & c_{1100} &   0   
\end{matrix} \right] \!\!
{\rm det} \! \left[ \begin{matrix}
 2 c_{1111} & c_{1110} & c_{1101} & c_{1011} & c_{0111} \\
   c_{1110} &   0      & c_{1100} & c_{1010} & c_{0110} \\
   c_{1101} & c_{1100} &    0     & c_{1001} & c_{0101} \\
   c_{1011} & c_{1010} & c_{1001} &     0    & c_{0011} \\
   c_{0111} & c_{0110} & c_{0101} & c_{0011} &   0
\end{matrix}\right].
\end{eqnarray*}
\end{smaller}

\noindent The squared factor is the $A$-discriminant of the octahedron,
and the $5 \times 5$-determinants are $A$-discriminants
of the two maximal cells in this subdivision.

\medskip
\noindent {\bf Facet 4:}
If $ w = - e_{0000} - e_{0001} $ then
${\rm in}_w(D_{2222})$ equals
\begin{smaller}
\begin{eqnarray*} &
(c_{0000} c_{1111}-c_{0001} c_{1110})^3  \\ &
(c_{1100} c_{1111}-c_{1101} c_{1110})
(c_{1010} c_{1111}-c_{1011} c_{1110})
(c_{0110} c_{1111}-c_{0111} c_{1110}) \\ &  \!\!\!
(c_{0000}^2 c_{1001} c_{1111} - c_{0000}^2 c_{1011} c_{1101}
- c_{0000} c_{0001} c_{1000} c_{1111} 
- c_{0000} c_{0001} c_{1001} c_{1110}  \\ &
+ c_{0000} c_{0001} c_{1010} c_{1101}
+ c_{0000} c_{0001} c_{1011} c_{1100}
+ c_{0001}^2 c_{1000} c_{1110} - c_{0001}^2 c_{1010} c_{1100})  \\ & \!\!\!
(c_{0000}^2 c_{0101} c_{1111}-c_{0000}^2 c_{0111} c_{1101}
-c_{0000} c_{0001} c_{0100} c_{1111}
-c_{0000} c_{0001} c_{0101} c_{1110} \\ &
+c_{0000} c_{0001} c_{0110} c_{1101}
+c_{0000} c_{0001} c_{0111} c_{1100}
+c_{0001}^2 c_{0100} c_{1110} - c_{0001}^2 c_{0110} c_{1100}) \\ & \!\!\!
(c_{0000}^2 c_{0011} c_{1111} - c_{0000}^2 c_{0111} c_{1011}
-c_{0000} c_{0001} c_{0010} c_{1111}
-c_{0000} c_{0001} c_{0011} c_{1110} \\ &
+c_{0000} c_{0001} c_{0110} c_{1011}
+c_{0000} c_{0001} c_{0111} c_{1010}
+c_{0001}^2 c_{0010} c_{1110}
-c_{0001}^2 c_{0110} c_{1010}).
\end{eqnarray*}
\end{smaller}
The quartic factors are the $A$-discriminants of the
three maximal cells.

\medskip
\noindent {\bf Facet 5:}
% $$ e_{0000}+e_{0001}+e_{0010}+e_{0100}+e_{0110}+e_{1000}+e_{1001} \geq 3 .$$
If $w = e_{0000}+e_{0001}+e_{0010}+e_{0100}+e_{0110}+e_{1000}+e_{1001} $ then
${\rm in}_w(D_{2222})$
equals 
$\, c_{0000} c_{0011} c_{0101} c_{1010} c_{1100}  \,$
times $\,(c_{0011} c_{1100} - c_{0101} c_{1010})^3 \,$
times
\begin{eqnarray*}
& c_{0011} c_{1100} c_{1111}-c_{0011} c_{1101} c_{1110}
-c_{0101} c_{1010} c_{1111} \\ & +c_{0101} c_{1011} c_{1110}
+c_{0111}c_{1010} c_{1101} - c_{0111} c_{1011} c_{1100}
\end{eqnarray*}
times the product of two larger factors,
each having $15$ terms of degree five
which are the $A$-discriminants of 
the two cells with $f$-vector $(10,23,21,8)$.

\medskip

In the remaining cases the tight spans
are three-dimensional (Figure 2), and
the initial form ${\rm in}_w(D_{2222})$ factors as a 
monomial times the product of the $A$-discriminants 
of the maximal cells in that  subdivision of the 4-cube.

\smallskip

\noindent {\bf Facet 6:}
%e_{0000}+e_{0001}+e_{0010}+e_{0011}+e_{0100}+e_{0110}+e_{1000}+e_{1001}\geq 4.
If $ w = e_{0000}+e_{0001}+e_{0010}+e_{0011}
+e_{0100}+e_{0110}+e_{1000}+e_{1001}$ then
${\rm in}_w(D_{2222})$ has $404$ terms and factors as a monomial
times the product of the $A$-discriminants of the four maximal cells:
\begin{eqnarray*}
& \,\,\,\,\,{\rm in}_w(D_{2222}) \quad = \quad  c_{0000} \cdot
c_{1100} \cdot
c_{1010}^4 \cdot
c_{0101}^4 \cdot 
\qquad \qquad \qquad \qquad \qquad 
\\ &
(c_{0101} c_{1010} c_{1111} - c_{0101} c_{1011} c_{1110}
-c_{0111} c_{1010} c_{1101} + c_{0111} c_{1011} c_{1100})
\\ &
(c_{0000} c_{0111} c_{1011} - c_{0001} c_{0111} c_{1010}
- c_{0010} c_{0101} c_{1011} + c_{0011} c_{0101} c_{1010})
\\ & \!\!\!\!
(c_{0000} c_{0101} c_{0111} c_{1110}
-c_{0000} c_{0111}^2 c_{1100}-c_{0010} c_{0101}^2 c_{1110} \\ & \qquad  \qquad
+c_{0010} c_{0101} c_{0111} c_{1100}
-c_{0100} c_{0101} c_{0111} c_{1010}
+c_{0101}^2 c_{0110} c_{1010} )
\\ & \!\!\!\!
(c_{0000} c_{1010} c_{1011} c_{1101}
-c_{0000} c_{1011}^2 c_{1100}
-c_{0001} c_{1010}^2 c_{1101} \\ & \qquad  \qquad
+c_{0001} c_{1010}c_{1011} c_{1100}
-c_{0101} c_{1000} c_{1010} c_{1011}
+c_{0101} c_{1001} c_{1010}^2).
\end{eqnarray*}

\smallskip
\noindent {\bf Facet 7:}
If $w = - e_{0000}-e_{0001}-e_{0010}-e_{0011}-e_{0100}-e_{1000}-e_{1100}$ then
\begin{eqnarray*}
& {\rm in}_w(D_{2222}) \quad = \quad
c_{1111} \cdot
c_{0011}^4  \cdot
c_{1100}^4 \cdot
\qquad \qquad \qquad \qquad \qquad \qquad \qquad 
\\ &
(c_{0000} c_{0011} c_{1100}
-c_{0001} c_{0010} c_{1100}
-c_{0011} c_{0100} c_{1000})
\\ &
(c_{0010} c_{1000} c_{1111}
-c_{0010} c_{1011} c_{1100}
-c_{0011} c_{1000} c_{1110}
+c_{0011} c_{1010} c_{1100})
\\ &
(c_{0010} c_{0100} c_{1111}
-c_{0010} c_{0111} c_{1100}
-c_{0011} c_{0100} c_{1110}
+c_{0011} c_{0110} c_{1100})
\\ &
(c_{0001} c_{1000} c_{1111}
-c_{0001} c_{1011} c_{1100}
-c_{0011} c_{1000} c_{1101}
+c_{0011} c_{1001} c_{1100})
\\ &
(c_{0001} c_{0100} c_{1111}
-c_{0001} c_{0111} c_{1100}
-c_{0011} c_{0100} c_{1101}
+c_{0011} c_{0101} c_{1100}).
\end{eqnarray*}

\smallskip
\noindent {\bf Facet 8:}
% $$2 \cdot e_{0000}+e_{0001}+e_{0010}+e_{0100}+ e_{1000} \leq 18. $$
If $w = - 2 \cdot e_{0000}-e_{0001}-e_{0010}-e_{0100}- e_{1000} $ then
${\rm in}_w(D_{2222})$ equals
\begin{eqnarray*}
c_{0000}^3  \cdot
c_{1111}^3 & \!\!\!\! \! \cdot \,
(c_{0000} c_{1100} c_{1111}
-c_{0000} c_{1101} c_{1110}-c_{0100} c_{1000} c_{1111}) \\ &
(c_{0000} c_{1010} c_{1111}
-c_{0000} c_{1011} c_{1110}-c_{0010} c_{1000} c_{1111}) \\ &
(c_{0000} c_{0110} c_{1111}
-c_{0000} c_{0111} c_{1110}-c_{0010} c_{0100} c_{1111}) \\ &
(c_{0000} c_{1001} c_{1111}
-c_{0000} c_{1011} c_{1101}-c_{0001} c_{1000} c_{1111}) \\ &
(c_{0000} c_{0101} c_{1111}
-c_{0000} c_{0111} c_{1101}-c_{0001} c_{0100} c_{1111}) \\&
(c_{0000} c_{0011} c_{1111}
-c_{0000} c_{0111} c_{1011}-c_{0001} c_{0010} c_{1111}) .
\end{eqnarray*}
\end{proposition}

These formulas show that $A$-discriminants for
subconfigurations of the cube appear
naturally as irreducible factors of 
leading forms of the hyperdeterminant.
This suggests the general problem of
studying $A$-discriminants and the related
tight spans  for various families of $0/1$-polytopes.
Other natural classes of configurations
for which a study of $A$-discriminants would be
interesting include generalized permutohedra,
Birkhoff polytopes, and reflexive polytopes.

\section{Vertices of the secondary polytope}

The $4$-cube is given by the $5 \times 16$-matrix
$A$ in the Introduction. For any generic vector
$w \in \RR^{16}$ we consider the tight span of $(A,w)$
whose cells are the bounded faces of the simple polyhedron
 $\,P_w=  \{ u \in \RR^5 : u \cdot A \leq w \}$. 
Each vertex of $P_w$ is indexed by five columns of $A$.
This collection of  $5$-tuples,
each regarded as a $4$-simplex, is
a {\em regular triangulation} $\Pi_w$ of the $4$-cube
(cf.~\cite[\S 5]{DRS}).

The regular triangulations of the 4-cube are in bijection
with the vertices of the secondary polytope.
As in Section 2, we define the secondary polytope  as
the Newton polytope $\N(E_{2222})$ of the {\em principal determinant} $E_{2222}$,
which is the following product of determinants associated to all
faces of the $4$-cube:
\begin{itemize}
\item $16$ linear factors $c_{ijkl}$ corresponding to the vertices,
\item  $24$ quadratic factors $D_{22}$  corresponding to the $2$-faces,
\item $8$ quartic factors $D_{222}$ corresponding to the facets,
\item the hyperdeterminant $D_{2222}$ corresponding to the solid $4$-cube.
\end{itemize}
Thus the principal determinant $E_{2222}$ is a polynomial
in the $16$ unknowns $c_{ijkl}$. of degree
$\,120 = 16 \cdot 1 + 24 \cdot 2 + 8 \cdot 4 + 1 \cdot 24 $.
The vertices of its Newton polytope $\N(E_{2222})$
are the GKZ vectors of the regular triangulations
of the 4-cube.

Recall from the first part of Theorem 2 that the 4-cube has
$87959448$ regular triangulations, or, equivalently,
$\N(E_{2222})$ has $87959448$ vertices.  Up to
the action of the Weyl group $B_4$, there are
$235277$ types of regular triangulations. Thus,
this is the number of combinatorially distinct
simple polyhedra $P_w$, and also
the number of generic tight spans.

\begin{figure}
\vskip-0.5in
%\vskip2in
\begin{tabular}{ccc}
\begin{tabular}{c}
\hskip-0.7in \includegraphics[scale=0.35]{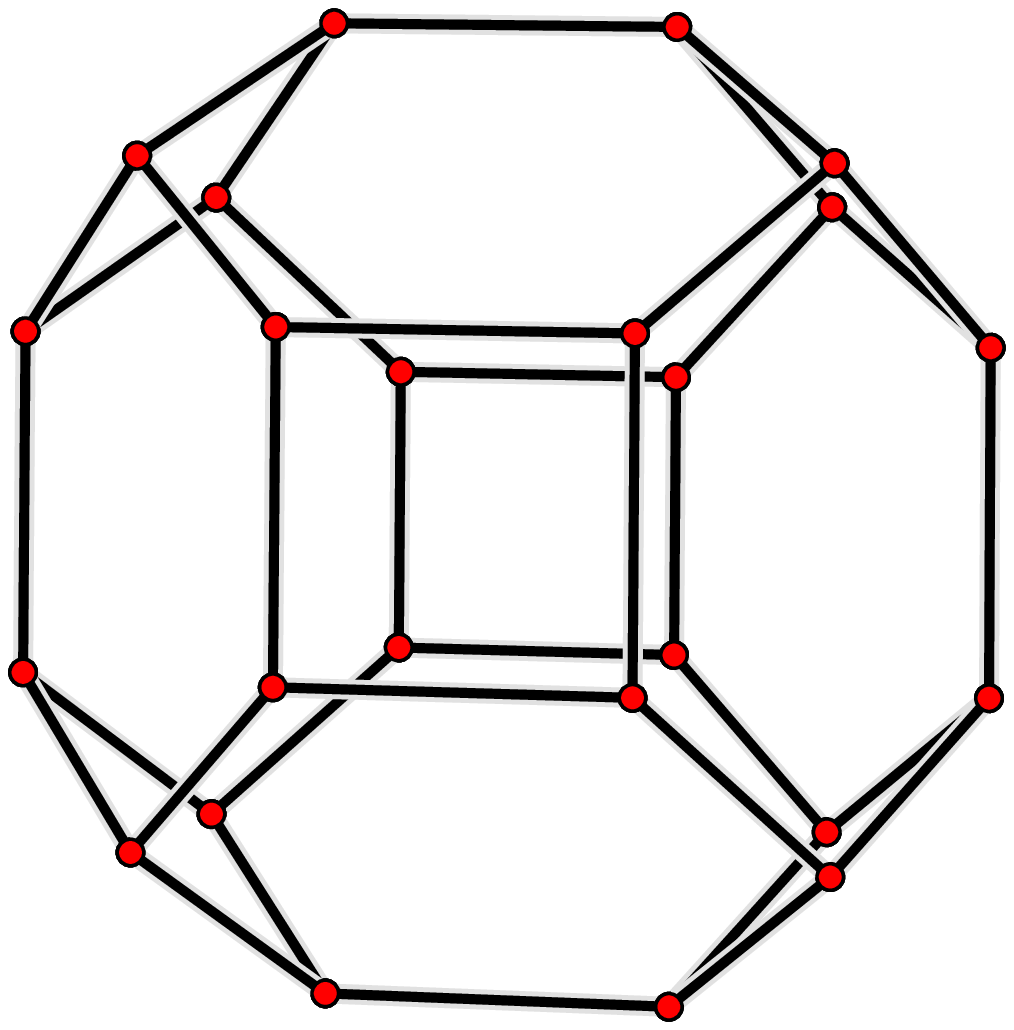}
\end{tabular}
& 
\begin{tabular}{c}
\hskip-1.6in \includegraphics[scale=0.45]{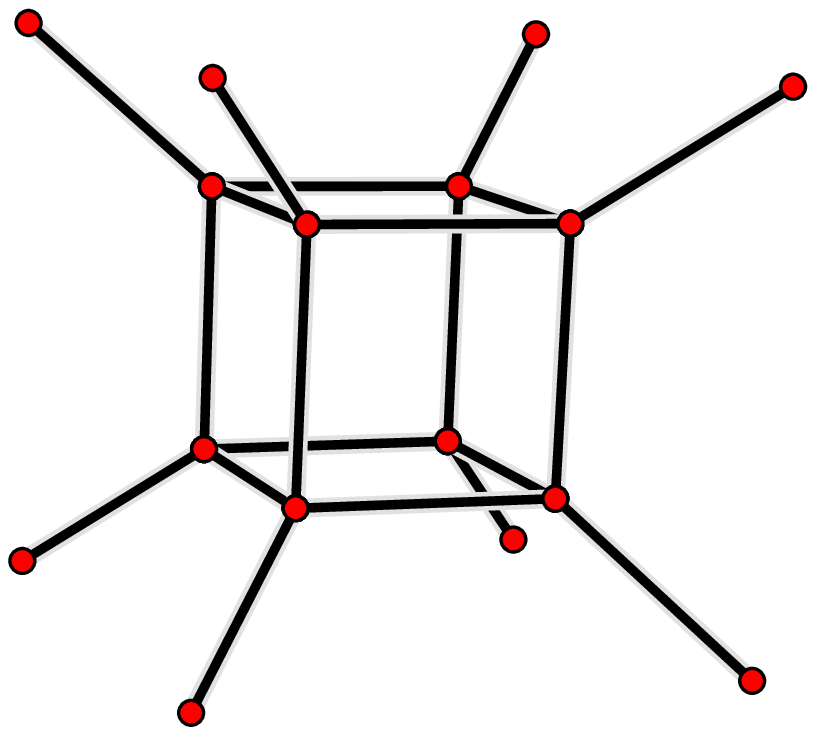}
\end{tabular}
& 
\begin{tabular}{c}
\hskip-1.8in \includegraphics[scale=0.4]{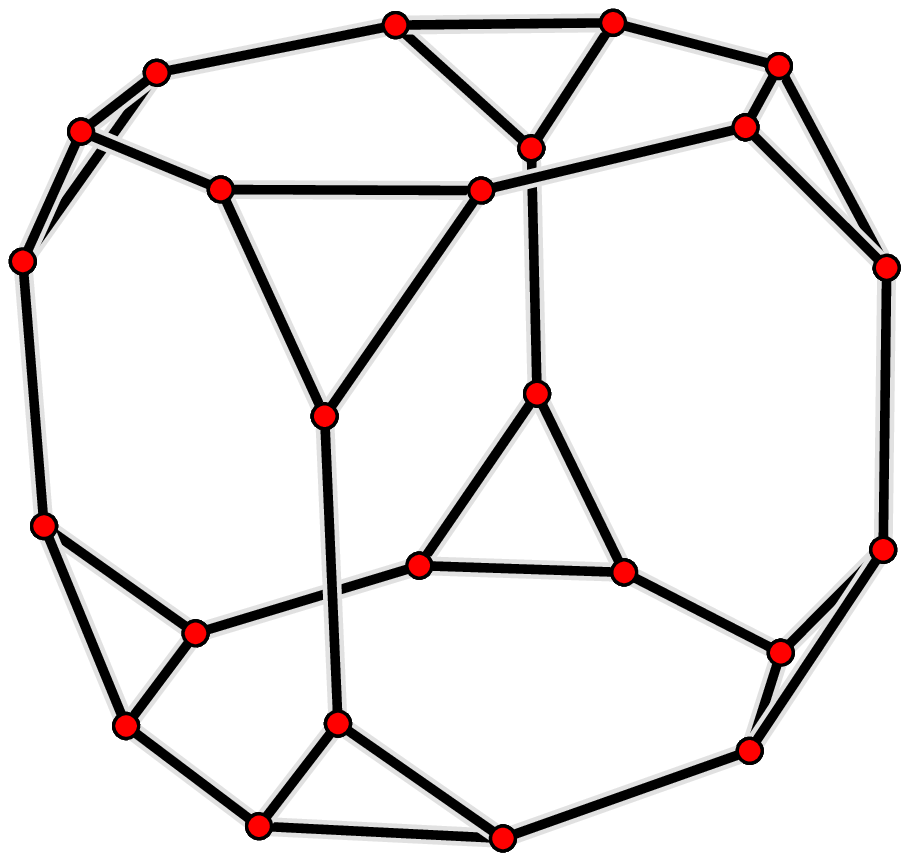}
\end{tabular}
\end{tabular}
\vskip-0.7in
\caption{Tight spans of the three triangulations
in Example \ref{mostsymmetric}}
\label{triangpic1}
\end{figure}

The main objective of this section is to provide
a detailed study of the regular triangulations of the
$4$-cube, with emphasis on their tight spans.
We present a range of results
which were computed from our data.

\begin{corollary}
The distribution of the $235277$ symmetry classes  of
vertices of the secondary polytope $\N(E_{2222})$ 
according to orbit size is as follows:
$$
\begin{array}{cccccccccccc}
8 & 16 & 24 & 32 & 48 & 64 & 96 & 192 & 384 \\
\underline{3} & 7 & 2 & 13 & 48 & 102 & 516 & 11357 & 223229 \\
\end{array}
$$
\end{corollary}

\begin{table}
\begin{tabular}{|c|ccc|c|}
\hline
\# maximal
& \multicolumn{3}{c|}{\# simplices of}
  &  {\begin{tabular}{c}\# $B_4$-orbits of\end{tabular}} \\
  simplices & volume 1 & volume 2 & volume 3 & triangulations\\
  \hline
  16&8&8&0&1\\ \hline
17&10&7&0&1\\ \hline
18&12&6&0&6\\
18&13&4&1&11\\ \hline
19&14&5&0&25\\
19&15&3&1&48\\ \hline
20&16&4&0&628\\
20&17&2&1&344\\ \hline
21&18&3&0&5847\\
21&19&1&1&1263\\ \hline
22&20&2&0&24499\\
22&21&0&1&1967\\ \hline
23&22&1&0&48648\\ \hline
24&24&0&0&151989\\
\hline
\end{tabular}
\medskip
\caption{Classification of regular triangulations
of the 4-cube according to the number and volumes of maximal simplices}
\label{tab:Nsimplices}
\end{table}

\begin{example}[\bf The most symmetric triangulations]
\label{mostsymmetric}
One of the triangulations with orbit size $8$
 is the {\em staircase triangulation}. Its
tight span is a solid permutohedron, and
its GKZ vector is $(4,6,6,24,6,4,4,$ $6,6,4,4,6,24,6,6,4)$.
Another is the unique triangulation that
uses 16 maximal simplices.
Eight simplices have volume one and
eight have volume two. Its GKZ vector is
  $(1,12,12,1,12,1,1,20,20,1,1,$ $12,1,12,12,1)$, and
its tight span is a $3$-cube with tentacles attached
to its vertices.  The third triangulation with orbit size
$8$ has a solid truncated $3$-cube as its tight span.
Its GKZ vector is $(3,8,8,3,8,3,3,24,24,3,3,8,3,8,8,3)$. 
The tight spans of these symmetric triangulations
are shown in
 Figure \ref{triangpic1}. \qed
\end{example}

Table 3 classifies the regular
 triangulations according to the number of maximal simplices used.
It is well-known that the minimal number is sixteen (cf.~\cite{Cot}).
  We see that any triangulation can only use at most one simplex of volume
three, since all such simplices contain the centroid of the cube 
as an interior point. There are $16$ simplices of volume three, 
all of which are $B_4$-equivalent; a representative is  $\,\{(1000),(1111),(0011),(0101), (0110)\}$.

\begin{figure}
\vskip-0.5in
\begin{tabular}{ccc}
\hskip-0.9in \includegraphics[scale=0.4]{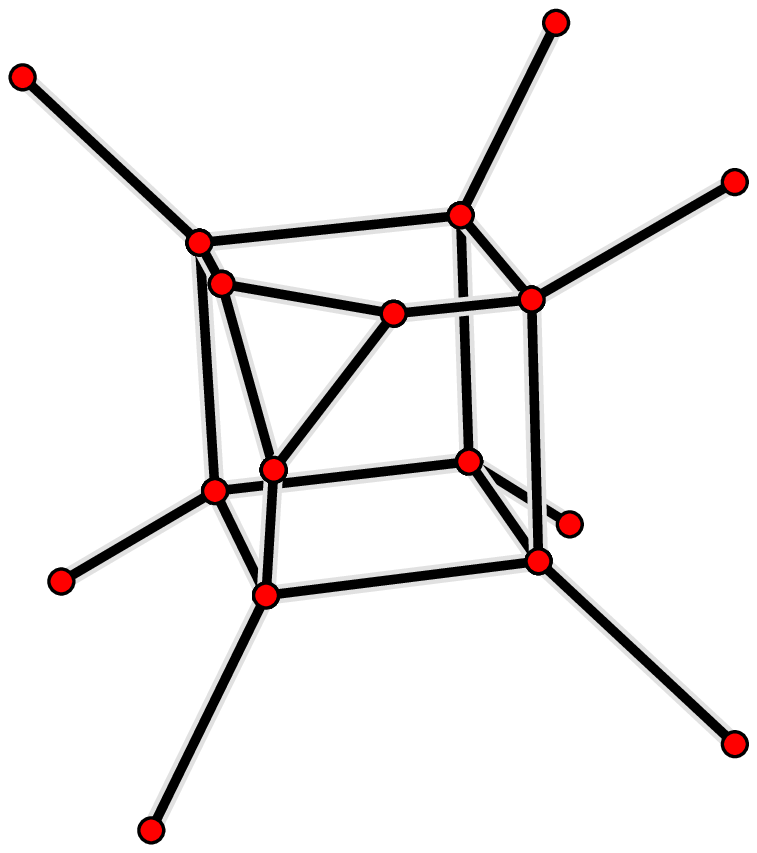}& 
\hskip-1.3in \includegraphics[scale=0.35]{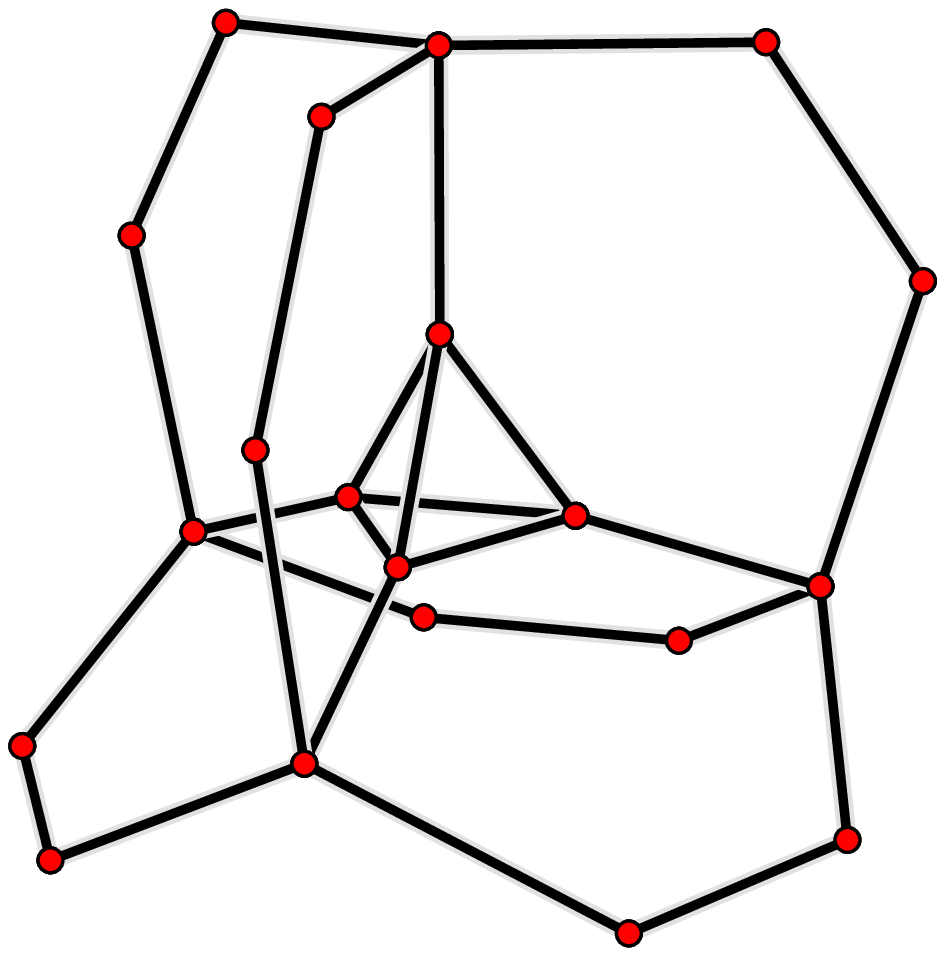}& 
\hskip-1in \includegraphics[scale=0.35]{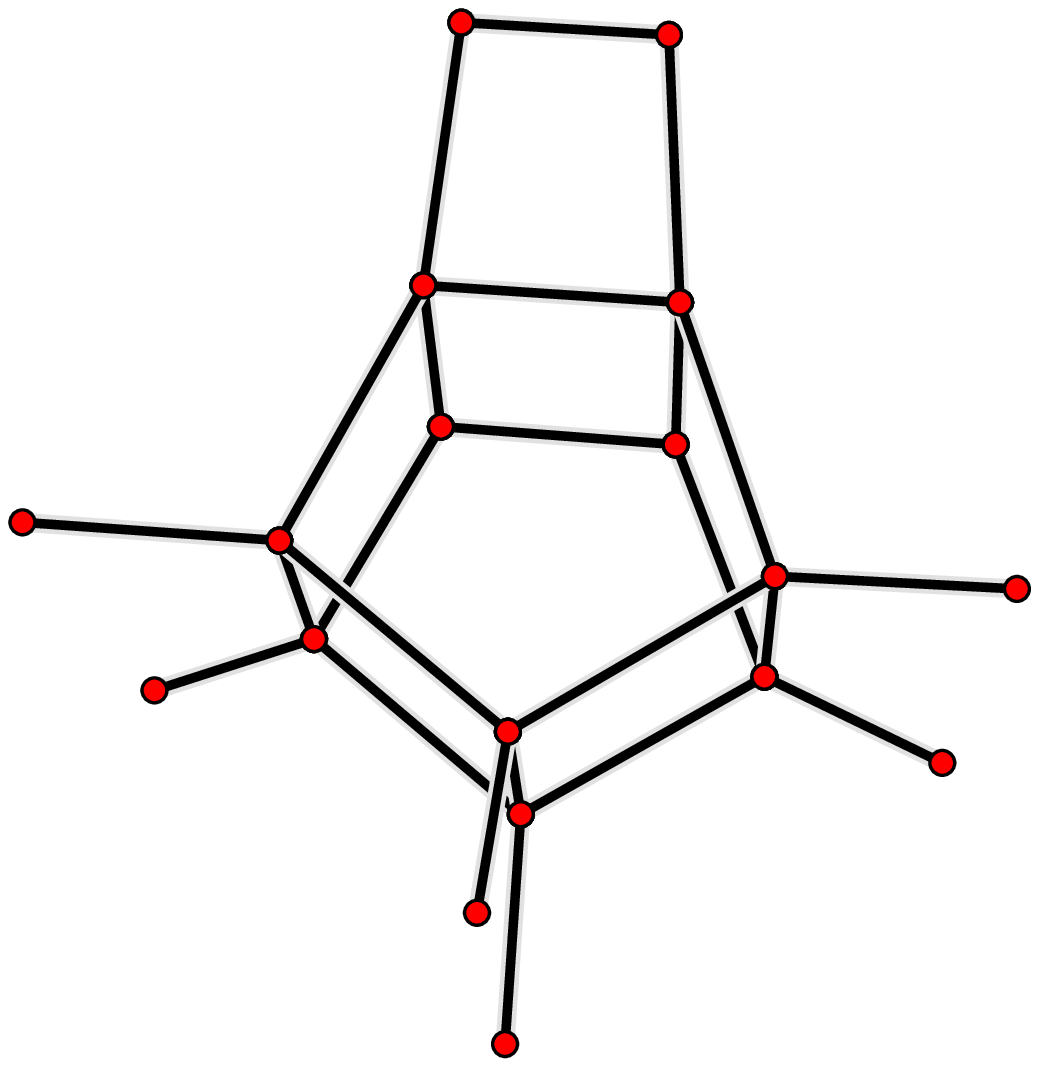}
\end{tabular}
\vskip-0.5in
\caption{The triangulations
in Examples \ref{sweet17} and \ref{smallesttight}}
\label{triangpic2}
\end{figure}

\begin{example}[\bf Triangulating the 4-cube with 17 simplices]
\label{sweet17}
Table 3 shows that there is a unique regular triangulation with
$17$ maximal simplices, ten of volume one and seven of volume two.
Its tight span, shown in
Figure \ref{triangpic2},
is a cube with a truncated vertex and tentacles
attached to the seven original vertices.
 Its GKZ vector is $(1,11,12,1,12,1,1,21,20,3,1,11,1,11,12,1)$. \qed
\end{example} 

\begin{table}
\begin{tabular}{|c|c|c|c|c|c|c|c|}
\hline
			 & \{3\} & \{3,2\} 	& \{3,1\} &\{3,2,1\} &\{2\} &\{2,1\}  & total\\
\hline
(16, 20, 6, 1)	 & 	&	& 1 & & & 					& 1 \\ 
\hline
(17, 22, 7, 1) 	&  & & 1 & & & &1 \\ 
\hline
(18, 24, 8, 1) 	& & & 4  &2 & & & 6 \\
(18, 23, 6, 0) 	&   &  &  &  &  &  11 & 11 \\ 
\hline
(19, 26, 9, 1) 	& & & 6 & 19 & &  & 25 \\
(19, 25, 7, 0) 	& & & & & & 48 & 48 \\
\hline
(20, 28, 10, 1) 	& & 1 & 23 & 209 & & &  233 \\
(20, 27, 8, 0) 	& & & & & 24 & 715 & 739 \\
\hline
(21, 30, 11, 1) 	& & 5 & 34 & 1372 & & & 1411 \\
(21, 29, 9, 0) 	& & & & & 392 & 5307 & 5699\\
\hline
(22, 32, 12, 1) 	& & 112 & 84 & 9342 & & &  9538\\
(22, 31, 10, 0) 	& & & & & 2156 & 14772 & 16928\\
\hline
(23, 34, 13, 1)	 & & 2116 & 100 & 46432 & & & 48648\\
\hline
(24, 36, 14, 1)	 & 125 & 27054 & & 124810 & & & 151989\\ 
\hline
total & 125 & 29288 & 253 & 182186 & 2572 & 20853 & 235277\\
\hline
\end{tabular}
\medskip
\caption{Distribution of the 235277 types of triangulations.  The rows correspond to the $f$-vectors and the columns correspond to the dimensions of maximal faces of the tight span.}
\label{tab:TSfvector}
\end{table}

Table \ref{tab:TSfvector} classifies the tight spans
according to their $f$-vectors and their {\em signature},
by which we mean the set of dimensions of the maximal cells.
Note that every triangulation of the 4-cube has at most one interior
edge (namely, a diagonal), so each tight span has
at most one $3$-dimensional cell.
The $f$-vector of a tight span is determined 
by the number of vertices and 3-cells,  i.e., 
the number of maximal simplices and diagonals used in the triangulation.   

\begin{example}[\bf The smallest tight spans of 
signature $\{3,2\}$ and $\{3,2,1\}$]
\label{smallesttight}
The unique tight span of signature $\{3,2\}$ with
$20$ vertices contains a solid tetrahedron with six hexagons attached, one to each of the tetrahedron's edges.  See Figure \ref{triangpic2}.  Its triangulation has 
GKZ vector $\,(5,2,2,15,2,15,15,4,24,5,5,2,5,2,$ $2,15)\,$
and lies in the 88th $D$-equivalence class. 

There are two tight spans of  signature
$\{3,2,1\}$ with $18$ vertices.  Their GKZ vectors are 
$( 1,9,13,2,12,1,1,21,19,6,1,9,1,11,12,1)$ and 
$(1,12,12,1,$ $11,1,1,21,19,1,1,13,5,10,10,1)$.
They are both in the largest $D$-equivalence class.
One of them is shown in Figure \ref{triangpic2}.  The other has a triangle and a tentacle in place of the rectangular maximal 2-face.
\qed
\end{example}

The first column in Table \ref{tab:TSfvector} is
especially interesting as it ties in 
with an active area of research in
combinatorial commutative algebra. Using the
language introduced in \cite{CHT}, the $d$-cube is a 
{\em Gorenstein polytope},
and the tight spans of signature $\{d-1\}$ are precisely its
{\em Gorenstein triangulations}.
The 3-cube has only one  Gorenstein
triangulation, corresponding to the last type in
Figure~\ref{cubetriangs}.

\begin{example}[\bf The Gorenstein triangulations of the 4-cube]
\label{gorenstein}
The 4-cube has precisely $125$ types of Gorenstein
triangulations. Two of them are pictured
in Figure \ref{triangpic1}.
 The Stanley-Reisner ideal of a Gorenstein 
triangulation corresponds to an initial ideal of the Segre variety
$\PP^1 \times\PP^1 \times \PP^1 \times \PP^1 \subset \PP^{15}$
that is squarefree and Gorenstein. The minimal free
resolution of the Alexander dual of this monomial ideal
is cellular and is supported on the tight span itself.
For tropical polytopes (Example \ref{ex1}) this was
shown in \cite{BY}, but the same holds
for unimodular regular triangulations of arbitrary
polytopes, including the 4-cube.
All $125$ Gorenstein tight spans
are {\em deformed permutohedra}, and it would
be interesting to study their combinatorics
in arbitrary dimensions $d \geq 4$.
\qed
\end{example}

\begin{table}
\begin{equation*}
\begin{array}{rrrrrrrrrrrrrrrr}
\multicolumn{4}{l}{87959448} &&&&&&&&&&&&\\
=& 32&\!\cdot\!\!& 2116 &+& 192&\!\cdot\!\!& 644 &+& 192&\!\cdot\!\!& 2752
&+&
64&\!\cdot\!\!& 30390\\
+& 384 &\!\cdot\!\!& 1550 &+& 192 &\!\cdot\!\!& 368 &+& 384 &\!\cdot\!\!&
3742 &+&
384 &\!\cdot\!\!& 716 \\
+& 192 &\!\cdot\!\!& 10384 &+& 192 &\!\cdot\!\!& 20552 &+& 96
&\!\cdot\!\!& 1444 &+&
384 &\!\cdot\!\!& 380 \\
+& 192 &\!\cdot\!\!& 4408 &+& 384 &\!\cdot\!\!& 1318 &+& 384 &\!\cdot\!\!&
940 &+&
192 &\!\cdot\!\!& 468 \\
+& 192 &\!\cdot\!\!& 2584 &+& 192 &\!\cdot\!\!& 18050 &+& 384
&\!\cdot\!\!& 4122 &+&
384 &\!\cdot\!\!& 11614 \\
+& 192 &\!\cdot\!\!& 9090 &+& 192 &\!\cdot\!\!& 5276 &+& 96 &\!\cdot\!\!&
2440 &+&
96 &\!\cdot\!\!& 22788 \\
+& 384 &\!\cdot\!\!& 3182 &+& 384 &\!\cdot\!\!& 1358 &+& 192 &\!\cdot\!\!&
1444 &+&
384 &\!\cdot\!\!& 760 \\
+& 192 &\!\cdot\!\!& 5140 &+& 384 &\!\cdot\!\!& 940 &+& 384 &\!\cdot\!\!&
4068 &+&
384 &\!\cdot\!\!& 2810 \\
+& 192 &\!\cdot\!\!& 7956 &+& 192 &\!\cdot\!\!& 15260 &+& 384
&\!\cdot\!\!& 1318 &+&
384 &\!\cdot\!\!& 3278\\
  +& 192 &\!\cdot\!\!& \underline{18964} &+& 96 &\!\cdot\!\!& \underline{18340} &+& 96
&\!\cdot\!\!& \underline{21192}
&+& 96 &\!\cdot\!\!& 196 \\
  +& 192 &\!\cdot\!\!& 392 &+& 384 &\!\cdot\!\!& 1616 &+& 192 &\!\cdot\!\!&
1428 &+&
192 &\!\cdot\!\!& 18042 \\
  +& 384 &\!\cdot\!\!& 454 &+& 384 &\!\cdot\!\!& 336 &+& 384 &\!\cdot\!\!&
112 &+&
384 &\!\cdot\!\!& 2112 \\
  +& 384 &\!\cdot\!\!& 780 &+& 192 &\!\cdot\!\!& 212 &+& 384 &\!\cdot\!\!&
336 &+&
192 &\!\cdot\!\!& 4624 \\
  +& 96 &\!\cdot\!\!& 100 &+& 192 &\!\cdot\!\!& 200 &+& 384 &\!\cdot\!\!&
1360 &+&
192 &\!\cdot\!\!& 338\\
  +& 384 &\!\cdot\!\!& 240 &+& 384 &\!\cdot\!\!& 288 &+& 384 &\!\cdot\!\!&
200 &+&
384 &\!\cdot\!\!& 96\\
   +& 192 &\!\cdot\!\!& 8770 &+& 384 &\!\cdot\!\!& 5852 &+& 384
&\!\cdot\!\!& 2180
&+& 384 &\!\cdot\!\!& 1692 \\
   +& 192 &\!\cdot\!\!& 1490 &+& 96 &\!\cdot\!\!& 392 &+& 48 &\!\cdot\!\!&
3344 &+&
384 &\!\cdot\!\!& 908 \\
   +& 192 &\!\cdot\!\!& 4056 &+& 384 &\!\cdot\!\!& 1392 &+& 96
&\!\cdot\!\!& 1152 &+&
192 &\!\cdot\!\!& 576 \\
   +& 384 &\!\cdot\!\!& 232 &+& 384 &\!\cdot\!\!& 480 &+& 96 &\!\cdot\!\!&
900 &+&
192 &\!\cdot\!\!& 1690 \\
   +& 384 &\!\cdot\!\!& 1096 &+& 384 &\!\cdot\!\!& 956 &+& 192
&\!\cdot\!\!& 1508 &+&
192 &\!\cdot\!\!& 1856 \\
   +& 384 &\!\cdot\!\!& 280 &+& 192 &\!\cdot\!\!& 676 &+& 96 &\!\cdot\!\!&
2304 &+&
384 &\!\cdot\!\!& 676 \\
   +& 192 &\!\cdot\!\!& 676 &+& 96 &\!\cdot\!\!& 2304 &+& 48 &\!\cdot\!\!&
1152 &+&
16 &\!\cdot\!\!& 141888 \\
   +& 192 &\!\cdot\!\!& 4600 &+& 96 &\!\cdot\!\!& 2888 &+& 192
&\!\cdot\!\!& 2700 &+&
64 &\!\cdot\!\!& 23584\\
   +& 192 &\!\cdot\!\!& 5140 &+& 192 &\!\cdot\!\!& 112 &+& 192
&\!\cdot\!\!& 424 &+&
64 &\!\cdot\!\!& 6570 \\
   +& 96 &\!\cdot\!\!& 400 &+& 192 &\!\cdot\!\!& 200 &+& 192 &\!\cdot\!\!&
580 &+& 96
&\!\cdot\!\!& 18340 \\
   +& 16 &\!\cdot\!\!& \underline{82832} &+& 48 &\!\cdot\!\!& 11764 &+& 48
&\!\cdot\!\!& 6200 &+& 96 &\!\cdot\!\!& 58150 \\
   +& 192 &\!\cdot\!\!& 4860 &+& 192 &\!\cdot\!\!& 96 &+& 192 &\!\cdot\!\!&
288 &+&
192 &\!\cdot\!\!& 352 \\
   +& 192 &\!\cdot\!\!& 1240 &+& 8 &\!\cdot\!\!& \underline{349555}&+& 32
&\!\cdot\!\!& \underline{64}&\\
\end{array}
\end{equation*}
\caption{Decomposition into $D$-equivalence classes}
\label{decompose}
\end{table}

We now return to our discussion of the hyperdeterminant $D_{2222}$.
Since $D_{2222}$ is a factor
of the principal determinant $E_{2222}$,
there is a natural many-to-one map from the vertices
of $\N(E_{2222})$ to the vertices of $\N(D_{2222})$.
Two regular triangulations are $D$-equivalent if they are
mapped to the same vertex of $\N(D_{2222})$.
For the 3-cube the $D$-equivalence classes
were described in Section~2.
Table \ref{decompose} shows the decomposition of 
the number $87959448$ of regular triangulations of the 
4-cube into the cardinalities of the orbits of $D$-equivalence classes.  
Thus the identity expressed in Table  \ref{decompose} is the 4-cube version
of the identity  (\ref{seventyfour}) for the 3-cube. The $111$ orbits of $D$-equivalence classes are listed in the order of the corresponding orbits of vertices of $\N(D_{2222})$ in Section 4.

In the triangulations of the 4-cube, the centroid $(1,\frac{1}{2},\frac{1}{2},\frac{1}{2},\frac{1}{2})$ is in the relative interior of one of the following three kinds of simplices: a 4-simplex of normalized volume 3, a tetrahedron of normalized volume 2, or a diagonal.  

\begin{example}[\bf The ``coefficient $-27$ class'']
\label{fatclass}
The triangulations which use a fixed simplex of normalized volume 
three form a single $D$-equivalence class
of cardinality $82832$.
This class corresponds to the vertex 
$ 0  0  0  4  0  4  4  0  8  0  0  0 0  0  0 4^{-27}_{16} $
of $\N(D_{2222})$. The corresponding
monomial in $D_{2222}$ has the largest absolute value of any coefficient,
namely $-27$, among all  the extreme monomials. \qed
\end{example}

\begin{example}[\bf The ``coefficient $\pm 16$ classes'']
\label{coeff16}
The triangulations in which the centroid is in a fixed tetrahedron of volume two form a $D$-equivalence class corresponding to one of the extreme monomials of $D_{2222}$ with coefficient $\pm 16$, which come in three symmetry classes underlined in Table \ref{decompose} and on page \pageref{Dvertices}.  There are 24 such tetrahedra, all lying in a single $B_4$ orbit.  The hyperplane spanned by the special tetrahedron contains eight vertices of the 4-cube, and there are four vertices on each side.  In each triangulation, one vertex from each side is joined to the special tetrahedron.  The Hamming distance between these two vertices can be two, three, or four, and they correspond precisely to the three $D$-equivalence classes with coefficient $\pm 16$.  The first underlined class corresponds to Hamming distance three where there are two choices on one side after the vertex on the other side is chosen.  Hence it has orbit size $24 \times 4 \times 2 = 192$.  The second and third underlined classes correspond to Hamming distances four and two respectively and have orbit size $24 \times 4 = 96$.  \qed
\end{example}

The other $D$-equivalence classes contain only triangulations using a diagonal of the 4-cube, and all those classes contain some unimodular triangulations.

\begin{figure}
\vskip-0.3in
\begin{tabular}{ccc}
\hskip-0.7in \includegraphics[scale=0.35]{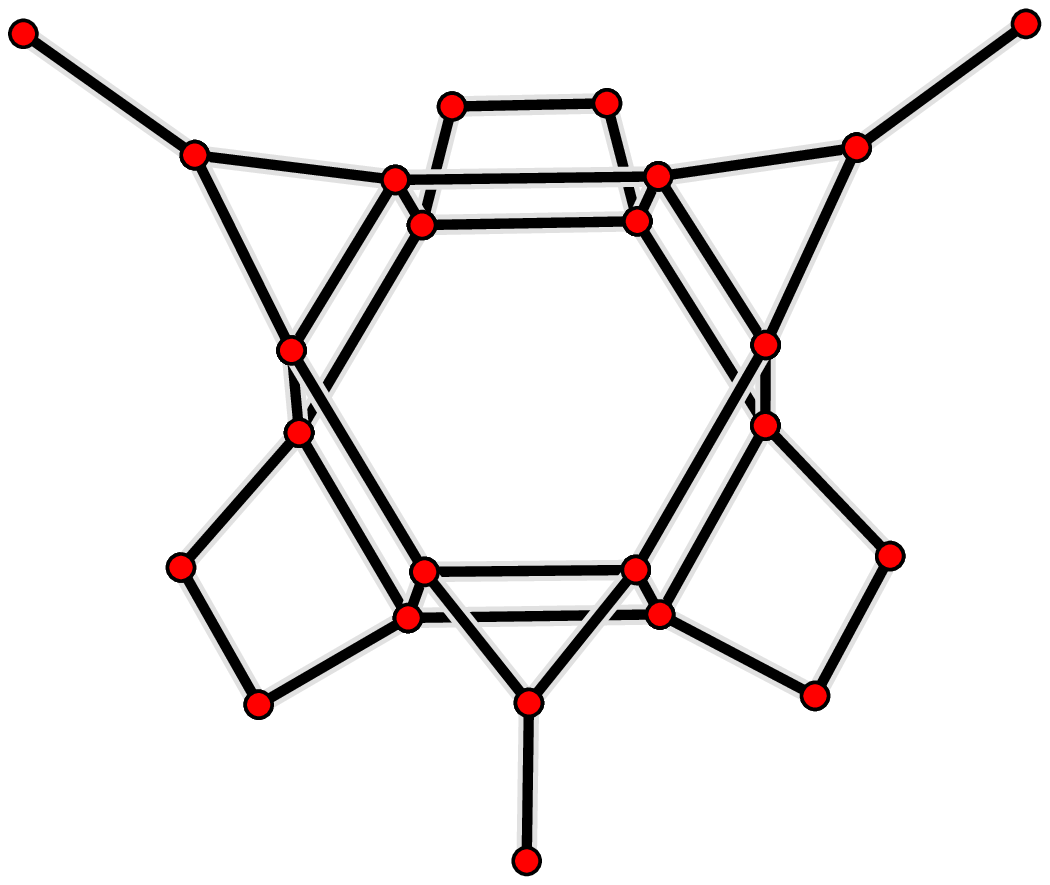}& 
\hskip-1.1in \includegraphics[scale=0.35]{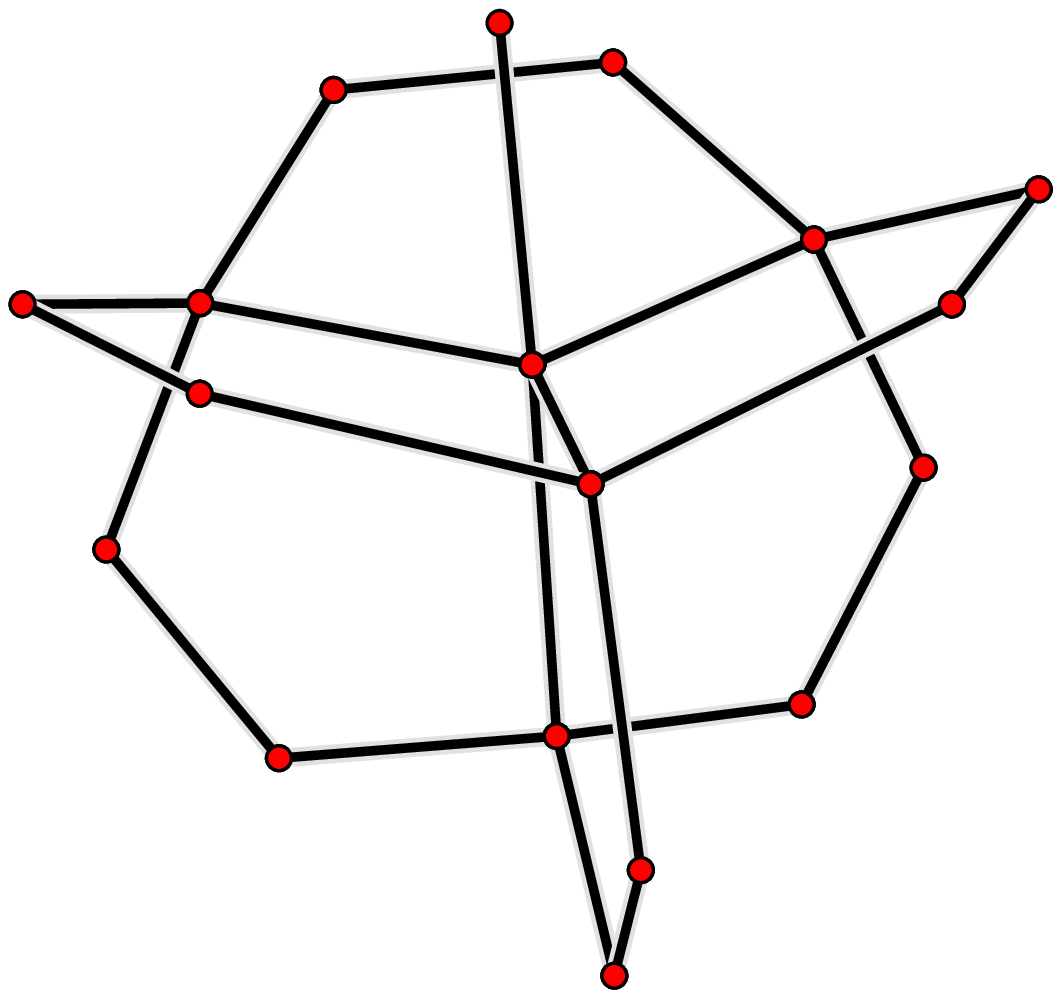}& 
\hskip-1.1in \includegraphics[scale=0.35]{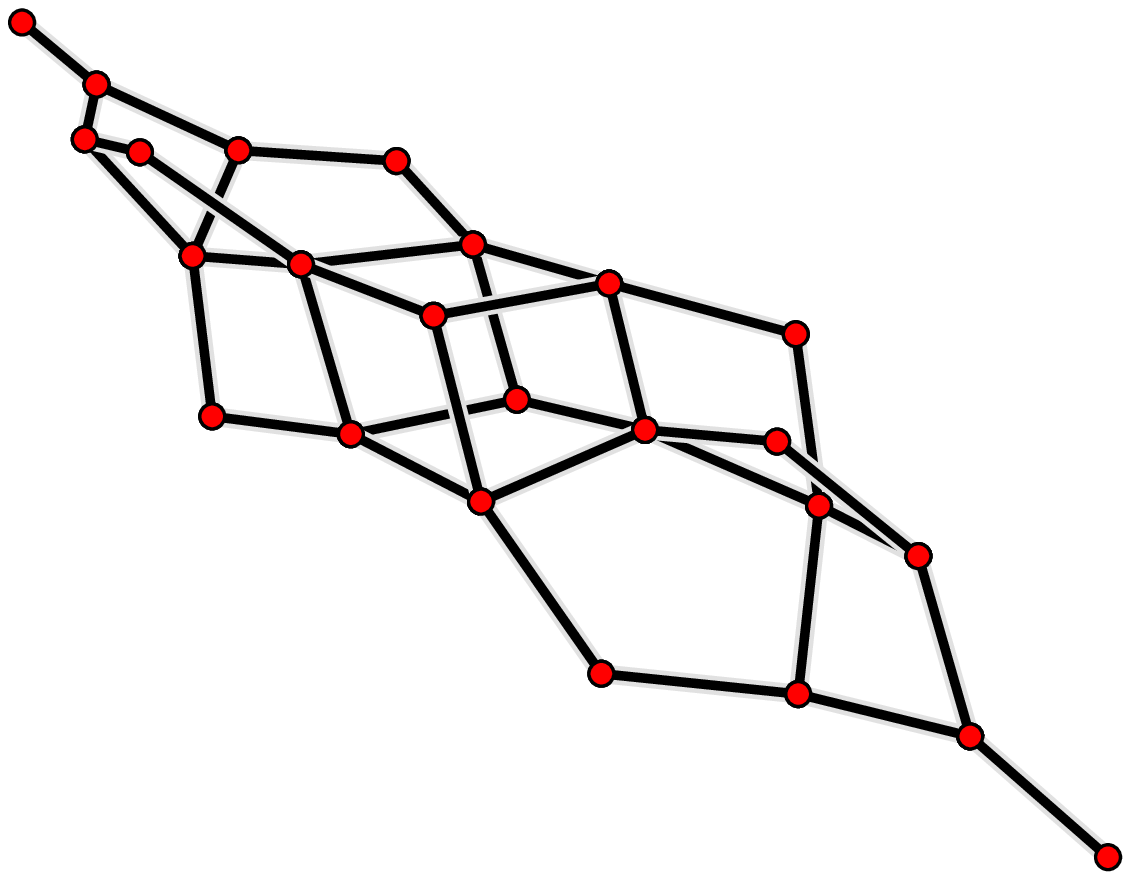}
\end{tabular}
\vskip-0.4in
\caption{Tight spans of triangulations
in Examples \ref{smallestclass}, \ref{fatclass} and \ref{delightful}}
\label{triangpic3}
\end{figure}

\begin{example}[\bf The smallest $D$-equivalence class]
\label{smallestclass}
This class has only $64$ triangulations and corresponds to the last summand
in Table \ref{decompose}.
The tight span of each of these 64 triangulations contains 
a hexagonal prism and six lower-dimensional pieces attached 
to alternating edges of the prism, as shown on
the left in Figure \ref{triangpic3}.
Each of these six lower-dimensional pieces
is either a rectangle or a triangle with a tentacle,
whence the number $\,2^6 = 64$. \qed
\end{example}

\begin{example}[\bf The largest $D$-equivalence class]
\label{largestclass}
This class has $349555$ triangulations and corresponds to
the vertex of $\N(D_{2222})$ which has the largest vertex figure and the smallest orbit size.
It contains all the triangulations in Examples
\ref{mostsymmetric}, \ref{sweet17} and \ref{gorenstein}.
\qed
\end{example}

  Using our data, it is easy to study other questions concerning regular
  triangulations of the 4-cube. Answering a question which
  was left open in \cite{SS}, we found that the
  4-cube does not possess any {\em delightful triangulations}.
  Such a triangulation would correspond to a square-free
  initial ideal of the first secant variety of
  $\PP^1 \times \PP^1 \times \PP^1 \times \PP^1 \subset \PP^{15}$,
  which has  dimension $9$ and degree $64$.
  For every regular triangulation $\Pi$ we computed the
  number of $10$-element subsets of the vertices
  of the 4-cube which are the union of vertices of two simplices
  in $\Pi$. That number is bounded above by $64$,
  and $\Pi$ would be delightful if equality holds.

\begin{example} [\bf The most delightful triangulation]
\label{delightful}
 We found that the maximum number  of
$10$-sets of vertices which are the unions of vertices of two simplices 
in a regular triangulation is $56$.
  This number is attained by two  triangulations which are both in
$D$-equivalence classes of type 11.  
Their GKZ vectors are $(1,4,4,11,6,13,9,12,12,9,13,6,11,4,4,1)$ and 
$(1,4,4,11,6,13, 9,$ $12, 12, 11,13, 4,9,4,6,1).$  Both of the corresponding tight spans have $f$-vector $(24,36,14,1)$ and consist of two maximal edges, eight maximal 2-faces, and one 3-cube.  The tight span of the first triangulation is shown in Figure \ref{triangpic3}. 
\qed
\end{example}  

\smallskip

\section{Facets of the secondary polytope}

The facets of the secondary polytope 
correspond to the coarsest regular polyhedral 
subdivisions $\Pi$.  
Here {\em coarsest} means that $\Pi$ refines no other proper subdivisions.  
We computed all such subdivisions for the $4$-cube.
This result gives an irredundant
 inequality representation of the secondary polytope $\N(E_{2222})$ of the 4-cube, analogous to that of the 3-cube in expression (\ref{NE222ineq}).  

Since the Newton polytope $\N(D_{2222})$ is a 
summand of the secondary  polytope $\N(E_{2222})$,
all facet normals of $\N(D_{2222})$ are also 
facet normals of $\N(E_{2222})$. Thus we already
know $268$ facets of the secondary polytope from
Section 4.  Furthermore the linear equations
$Ax = (120,60,60,60,60)^T $ hold on $\N(E_{2222})$. 

In order to compute the facets of the secondary polytope of the 4-cube, we first wrote a program to compute the GKZ vector of each regular triangulation of the 4-cube. The information about the triangulations themselves was obtained using TOPCOM.  We then used the following method to compute the facets of the secondary polytope, which can be applied to any point configuration with symmetry.

Since $\N(E_{2222})$ is an $11$-dimensional polytope with $87959448$
vertices, computing its facets is currently out of scope of any
general polyhedral software.  Thus, to compute the facets, we took advantage of the relationship between the secondary polytope and subdivisions of the 4-cube. Namely, we
know that each vertex of $\N(E_{2222})$ corresponds to a regular triangulation of the 4-cube,
and, by \cite[\S 7.2.C]{GKZ},  each edge $(v,w)$ of $\N(E_{2222})$ 
corresponds to the circuit connecting the two triangulations represented by the vertices $v$ and $w$.

Note that, given any triangulation of the 4-cube, there are relatively few 
circuits which contain a face of the triangulation.  As a result,
each tangent cone of $\N(E_{2222})$ is generated by relatively few vectors.
Hence we computed the facets of $\N(E_{2222})$ by equivalently computing the 
facets of its tangent cones.  Although there are $87959448$ tangent cones, they come in $235227$ orbits under the $B_4$-action.  Thus, we only
computed the facets of $235227$ tangent cones; all other facets of 
$\N(E_{2222})$ were obtained by applying the $B_4$-action.

We found that $\N(E_{2222})$ has $80876$ facets in $334$ orbits.  
Thus the 4-cube admits exactly $80876$ coarsest regular subdivisions. 
The distribution of the types of coarsest regular subdivisions according to orbit size is
$$
\begin{array}{cccccccccc}
8 & 12 & 16 & 24 & 32 & 48 & 64 & 96 & 192 & 384 \\
2 & 1 & 4 & 2 & 4 & 14 & 16 & 26 & 132 & 133 \\
\end{array}
$$

\begin{example}[\bf The most symmetric coarsest regular subdivisions]
The three most symmetric coarsest regular subdivisions come from 
facets of $\N(D_{2222})$.  
The two coarsest subdivisions with orbit size 8 
correspond to Facet 3 and Facet 8
in Sections 4 and 5. Their tight spans are a 
segment and an octahedron respectively.  The unique subdivision with orbit size 12 corresponds to Facet 2 in Sections 4 and 5.
Its tight span is also a line segment.
\qed
\end{example}

The distribution of the orbits according to the number of maximal cells in the subdivision and the dimensions of maximal faces of the  tight span are shown in Table \ref{NEfacets}.  
A coarsest subdivision of the 4-cube can have up to 13 maximal cells.  Here it is no longer the case that the number of vertices and 3-faces determine the $f$-vector of the tight span.  There are 33 types of $f$-vectors. 

\begin{table}
\begin{tabular}{|c|c|c|c|c|c|}
\hline
 & \{3\} & \{3,2\} & \{2\} & \{1\} & total \\
\hline
2& & & &4&  4\\
3& & &5& & 5\\
4&3& &7& & 10\\
5&4&3&17& & 24\\
6&5&14&12& & 31\\
7&13&21&39& & 73\\
8&7&31&51& & 89\\
9&7&26&24& & 57\\
10&2&14&21& & 37\\
11&1&1& & & 2\\
12&1& & & & 1 \\
13&1& & & & 1 \\
\hline
total & 44 & 176 & 110 & 4 & 334\\
\hline
\end{tabular}
\vskip .2cm
\caption{Distribution of the 334 types of coarsest regular subdivisions.  
The rows correspond to the number of maximal cells 
 and the columns correspond to the 
signatures of the tight span.}
\label{NEfacets}
\end{table}

\begin{figure}
\vskip-0.5in
\begin{tabular}{ccc}
\begin{tabular}{c}
\includegraphics[scale=0.55]{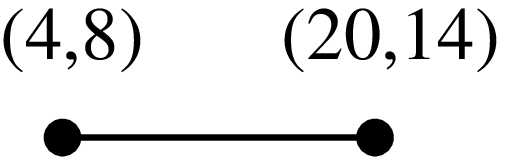}
\end{tabular}
&
\hskip-0.2in
\begin{tabular}{c}
\includegraphics[scale=0.3]{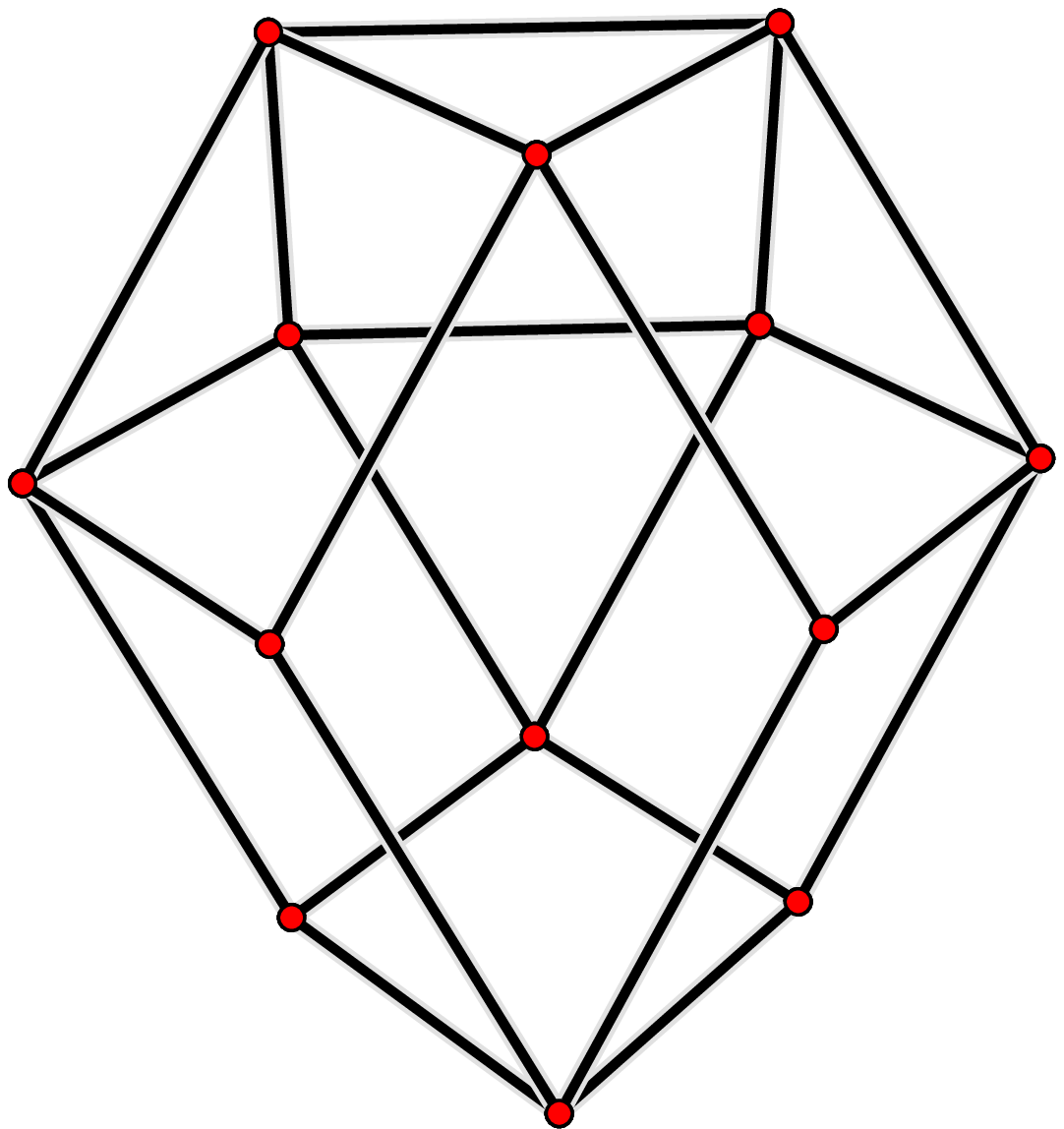}
\end{tabular}
&
\hskip-0.8in
\begin{tabular}{c}
\includegraphics[scale=0.35]{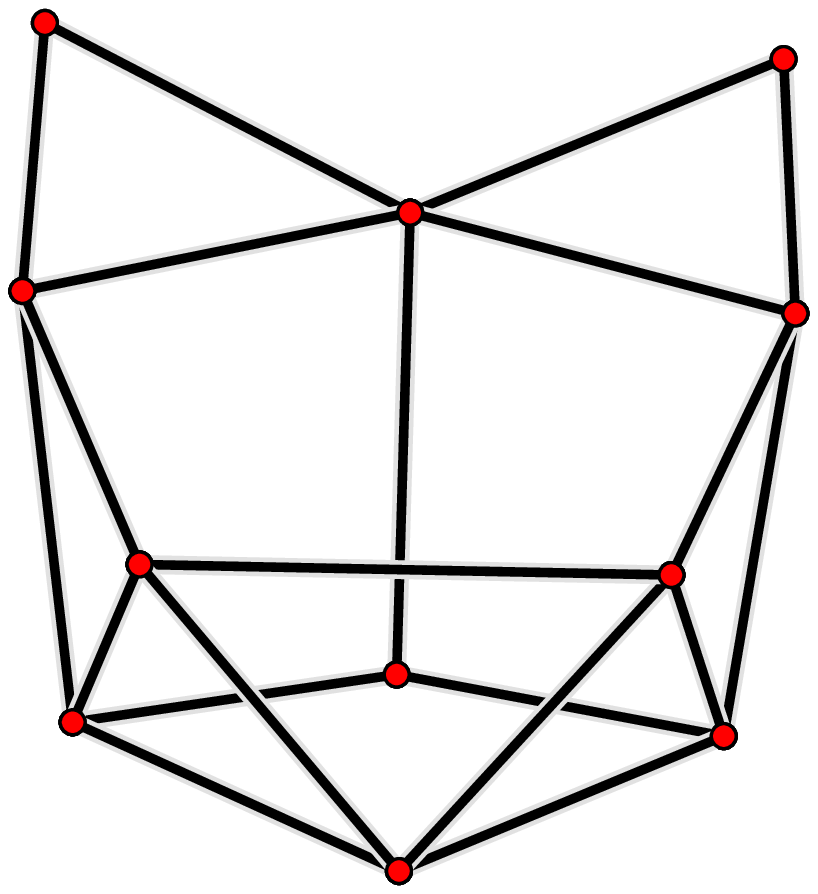} \\
\end{tabular}
   \end{tabular}
   \vskip-0.3in
\caption{The tight spans for  Examples \ref 
{split4}, \ref{largestcoarse}, and \ref{nonpurecoarse}}
\label{coarsefig}
\end{figure}

\begin{example}[\bf The missing split]
\label{split4}
The 4-cube admits 4 splits, i.e., subdivisions obtained by slicing 
the 4-cube with a hyperplane, whose tight spans are line segments.  Three of them correspond to Facets 1, 2, and 3 of $\N(D_{2222})$.  The remaining split subdivides the 4-cube into 2 cells with 8 and 14 vertices each, having normalized volumes 4 and 20 respectively.  The corresponding facets of $\N(E_{2222})$ are given by inequalities of the form 
$\,
x_{0000} + x_{0001} \geq 5 $.
\qed
\end{example}

\begin{example}[\bf The largest tight span]
\label{largestcoarse}
There is a unique coarsest regular subdivision with 13 maximal cells.  Its tight span is a  3-dimensional polytope with $f$-vector $(13,24,13,1)$, as shown in Figure \ref{coarsefig}.  The corresponding facets are given by inequalities 
like $\, x_{0000}+ x_{1111}  - x_{0011} \leq 44 $.
 \qed
\end{example}

\begin{example}[\bf The largest tight span with non-pure dimension]
\label{nonpurecoarse}
Almost half of the coarsest regular subdivisions have a pure dimensional tight span.  The largest non-pure dimensional tight span has $f$-vector $(11,20,11,1)$.  It contains one 3-dimensional face and two maximal 2-faces.  A figure is depicted in Figure \ref{coarsefig}.
The corresponding facets are given by inequalities like
$$
x_{0000}+x_{1010}+x_{1111} - x_{1000} - x_{1001} - x_{1011} \leq 47 .
$$
\vskip -.5cm
\qed
\end{example}

All the 334 facet inequalities and all the data we have discussed in this paper are available on our website
\url{bio.math.berkeley.edu/4cube/}.

\bigskip
\bigskip

\noindent {\bf Acknowledgements:}
This work began during Debbie Yuster's
Fall 2005 visit to the UC Berkeley Mathematics Department.
Peter Huggins was supported by an ARCS Foundation Fellowship.
Bernd Sturmfels was partially supported by the National Science
Foundation (DMS-0456960), and
Josephine Yu was supported by an NSF
Graduate Research Fellowship.  We also thank Francisco Santos for helpful discussions and  comments.

\bigskip

\noindent {\bf Authors' addresses:}

\bigskip

\noindent  Peter Huggins, Bernd Sturmfels and Josephine Yu,

  Department of Mathematics, University of California,  

   Berkeley, CA 94720, USA,

{\tt [phuggins,bernd,jyu]@math.berkeley.edu}

\medskip

\noindent Debbie S. Yuster, 

Department of Mathematics, Columbia University, 

New York, NY 10027, USA,

{\tt debbie@math.columbia.edu}

\end{document}